\documentclass[]{article}
\usepackage{geometry}
\usepackage{amsmath,amssymb,amsfonts,amsthm,amsopn,dsfont, mathrsfs}
\usepackage{bm}
\usepackage{graphicx}
\usepackage{float}
\usepackage[backend=bibtex,maxnames=5]{biblatex}

\newtheorem{rem}{Remark}

\usepackage{standalone}
\usepackage{tikz}
\usepackage{pgfplots}
\usepackage{pgfplotstable}
\pgfplotsset{compat=newest}
\usepackage{subcaption}
\usepackage{colortbl}
\usepackage{booktabs}

\usepackage[normalem]{ulem}

\usepackage[T1]{fontenc}

\usepackage{xcolor}
\usepackage{stmaryrd}

\usepackage[ruled,vlined,linesnumbered]{algorithm2e}
\usepackage{multirow}
\newcommand{\xx}{\bm{x}}

\newcommand{\I}{\mathcal{I}}

\newcommand{\tr}[2]{{#1}_{|{#2}}}

\newcommand{\phic}{\phi_c}
\newcommand{\phil}{\phi_l}
\newcommand{\phimax}{\phi_{max}}
\newcommand{\phiz}{\phi_{0}}
\newcommand{\phiecm}{\phi_m}

\newcommand{\D}{\mathcal{D}}

\title{Modeling tumor growth with variable mass and angiogenesis-driven perfusion through a 3D-1D coupled framework}
\author{Chiara Giverso \thanks{Dipartimento di Scienze Matematiche G. L. Lagrange, Politecnico di Torino, Italy (chiara.giverso@polito.it), Member of GNFM INdAM Group} \and Denise Grappein\thanks{MOX, Department of Mathematics, Politecnico di Milano, Italy, Member of GNCS INdAM Group ({denise.grappein@polimi.it})} \and  Stefano Scial\`o, \thanks{Dipartimento di Scienze Matematiche G. L. Lagrange, Politecnico di Torino, Italy (stefano.scialo@polito.it), Member of GNCS INdAM Group} }

\begin{document}
	\maketitle
	\begin{abstract}
    {Tumor growth beyond a critical size relies on the development of a functional vascular network, which ensures adequate oxygen and nutrient supply. In this work, we present a modeling framework based on an optimization-based 3D–1D coupling strategy to simulate perfusion in a tumoral tissue with growing mass, interacting with a dynamically evolving capillary network. The tumor is described as a multiphase system including tumor cells and interstitial fluid, governed by a nonlinear PDE system for cell volume fraction, pressure, oxygen, and VEGF, and discretized via finite elements. Capillary growth is tackled using a continuous-discrete hybrid tip-tracking approach. The vascular geometry is updated over time according to angiogenic signals, and coupled to the tissue model through a constrained optimization formulation that enforces fluid and nutrient exchange via interface variables. A sensitivity analysis using the Morris elementary effect method identifies key parameters influencing system behavior. Results highlight the critical role of vascular development in regulating tissue perfusion and tumor progression.
    Overall, the proposed numerical approach provides a versatile tool for investigating tumor–vascular interactions and can support further quantitative analysis of angiogenesis and tumor perfusion dynamics.\\
    
    \noindent \textbf{Keywords:} Optimization-based 3D-1D coupling, vascularized tumor growth modeling, mathematical model of angiogenesis, modelling tumor perfusion\\
    \textbf{MSC[2020]}   65N30, 35Q92, 92B05, 92C17 }
	\end{abstract}

\section{Introduction} 
Tumor-induced angiogenesis is a key mechanism supporting cancer progression, as it enables the delivery of oxygen and nutrients to proliferating malignant cells. Indeed, from a biological viewpoint, tumor angiogenesis is triggered by hypoxia. As the tumor mass expands beyond the critical size of $1-2 \rm mm$, the local oxygen concentration decreases due to insufficient perfusion and chemical diffusion \cite{Folkman_1985, Folkman_1987}. Hypoxic tumor cells respond by upregulating the production of pro-angiogenic factors, most notably the Vascular Endothelial Growth Factor (VEGF). The resulting VEGF gradients activate nearby endothelial cells (ECs), composing the capillary walls. Activated ECs degrade the basement membrane, migrate toward the tumor, proliferate, and eventually form new capillary sprouts. 
These sprouts elongate through the coordinated migration of specialized endothelial cells, called tip cells, which lead the advancing front, and the proliferation of trailing stalk cells, which form the vessel lumen. During their progression, sprouts may branch to generate new vascular segments and can eventually anastomose with neighboring vessels.
The whole process emerges from a complex interplay between vascular dynamics, biochemical signaling, and the mechanical and metabolic activity of the surrounding tissue \cite{Folkman_1985, Folkman_1995}. Understanding such multiscale interactions is crucial for improving the predictive capability of in silico models and for supporting the development of targeted therapies.

Over the past decades, a large variety of mathematical models  of angiogenesis has been proposed \cite{SciannaPreziosi, Gomez_review}. Early continuum approaches described endothelial cell density and chemoattractant fields through reaction-diffusion and chemotaxis equations (see, e.g., \cite{ByrneChaplain, Deakin, chaplain}). While these models successfully capture averaged vascular densities, they do not explicitly reproduce the geometry of the developing capillary network. To overcome this limitation, hybrid discrete--continuous models have been introduced, where tip cells are tracked individually while chemical fields are treated as continuous variables \cite{andersonchaplain, chaplain, Plank, Dougall2, sun, Milde}. These approaches provide a more realistic description of vascular morphology, although the coupling with interstitial flow and mass transport is often simplified. 
Alternative approaches based either on phase-field formulations (see, e.g., \cite{Gomez, Gomez2017, Travasso, Xu_Gomez}) or on sharp-interface descriptions \cite{GiversoCiarletta} have been proposed to still capture the evolving vessel morphology while preserving a fully continuous modelling framework.

More recently, increasing attention has been devoted to models that explicitly incorporate blood flow and transvascular exchanges \cite{BGGPS2023, Rieger2006, Rieger2008, Rieger2013, Fritz}. 
Many of these approaches build upon well-established models originally developed for haemodynamics in fixed microvascular networks and later adapted to angiogenic settings.

In this context, three-dimensional/one-dimensional (3D-1D) coupled formulations have proven particularly effective for describing blood flow and oxygen transport in complex vascular geometries (see, e.g., \cite{Koppl, cattzun0, cattzun}), providing a suitable compromise between physiological accuracy and computational efficiency. 
These approaches have also been successfully extended to angiogenesis and tumor growth models, allowing the vascular network to evolve while maintaining a consistent representation of intravascular flow and transvascular exchange (see, e.g., \cite{Fritz, BGGPS2023}).
Despite these advances, the consistent integration of detailed blood flow modeling with multiphase descriptions of the surrounding tumor tissue, as well as hypoxia-driven VEGF production, remains only partially explored.

In particular, the model proposed in \cite{BGGPS2023} combines hybrid tip-tracking angiogenesis with optimization based fluid and chemical exchange between a growing capillary network and the surrounding medium. However, in that setting the interstitial tissue was assumed to have constant density and the active growth of the tumor mass was not included.
From the biological standpoint, this limitation is significant. Tumor cells continuously proliferate, consume oxygen, and modulate VEGF production in response to the local microenvironment. These processes generate a strong feedback between tissue evolution and vascular remodeling that is known to play a central role in tumor progression. Although several multiphase models of tumor growth have been proposed in the literature (see, e.g., \cite{Preziosi, Wise2008, Frieboes2010}), their consistent integration with discrete angiogenesis and detailed vessel-tissue exchange mechanisms is still limited. A recent attempt can be found in \cite{Vavourakis}; however, a fully rigorous coupling between the tissue dynamics and the 3D-1D reduction of the vascular flow is still lacking.

The present work aims to address this issue by extending the framework introduced in \cite{BGGPS2023}. Specifically, we incorporate the evolution of the tissue surrounding the vessels, which in the considered application represents a growing tumor treated as a multiphase porous medium. The model accounts for tumor cell proliferation, oxygen transport and consumption, and endogenous VEGF production regulated by hypoxia. This results in a two-way coupling between the vascular network and the tumor microenvironment: on the one hand, the vascular network supplies oxygen and nutrients to the tissue; on the other hand, the evolving tumor microenvironment regulates VEGF production and thus guides angiogenic sprouting.

From the modeling viewpoint, the proposed framework therefore combines (i) hybrid discrete tip-tracking angiogenesis, (ii) 3D-1D fluid and nutrient exchange, and (iii) multiphase tumor tissue dynamics within a unified computational setting. This integrated approach allows us to investigate how vascular remodeling and tumor growth mutually influence each other.

The remainder of the paper is organized as follows. After the description of the used notation and of the main modeling assumptions in Section~\ref{sec:Notation}, Section~\ref{sec:model} introduces the multiphase tissue model and its coupling with the vascular network, while Section~\ref{sec:numerics} presents the numerical strategy. Section~\ref{sec:results} reports a systematic sensitivity analysis with respect to the model parameters together with representative numerical experiments highlighting the impact of the proposed feedback mechanisms. Conclusions and perspectives are finally drawn in Section~\ref{sec:conclusions}.

	\section{Notation and main assumptions}\label{sec:Notation}
In the following, we introduce the notation employed in this paper, following the same conventions adopted in previous work \cite{BGGPS2023} to ensure clarity and consistency.
	
Let $\Omega\subset\mathbb{R}^3$ be a tissue sample with boundary $\partial \Omega$. We consider a time interval $[0,T]$ partitioned in $K$ sub-intervals $\I_k=(t_{k-1},t_k]$, $k=1,...,K$. At time $t_0=0$ we consider an initial capillary network $\Sigma^0\subset \Omega\subset \mathbb{R}^3$, and we denote by $\D^0=\Omega \setminus \overline{\Sigma^0}$ the interstitial volume surrounding the capillaries. For the sake of simplicity, we allow for vessel growth but not for regression. Namely, denoting by $\Sigma^k=\Sigma(t_k)$ the capillary network at time $t_k$, we assume that $\Sigma^{k-1}\subseteq \Sigma^k$ for any $k=1,...,K$. The interstitial volume at time $t_k$ is defined as $\D^k=\Omega \setminus \overline{\Sigma^k}$.

We assume the network $\Sigma^k$ to be composed by thin tubular vessels of constant radius $R\ll~\mathrm{diam}(\Omega)$, the generalization to a radius varying along the network being straightforward.
	The boundary of $\Sigma^k$ is denoted by $\partial \Sigma^k$ and is split into four parts, i.e. $$\partial \Sigma^k= \Gamma^k\cup\left( \partial\Sigma_{in}\cup\partial\Sigma_{out}\cup\partial\Sigma_d^k\right).$$ The symbol $\Gamma^k$ is used to denote the lateral surface of $\Sigma^k$, whereas $\partial \Sigma_{in}, \partial \Sigma_{out}$ and $\partial \Sigma_d^k$ refer to the end cross sections, which we assume to lie either completely in $\Omega$ or completely on $\partial \Omega$. The set $\partial \Sigma_d^k$ collects the end sections lying in the interior of $\Omega$, whereas $\partial\Sigma_{in}\subset \partial \Omega$ and $\partial\Sigma_{out}\subset \partial \Omega$ are the union, respectively, of the inlet and outlet cross-sections, the nomenclature referring to blood velocity.
        The sets $\partial \Sigma_{in}$ and $\partial \Sigma_{out}$ are assumed to be fixed in time, excluding the possibility for any growing vessel to reach the boundary of the domain. In practice, $\Omega$ bounds the entire vessel network generated by angiogenesis at any time, encompassing the region in which the phenomena of interest take place.  This assumption is introduced mainly to simplify the notation so that each capillary tip in the growing network can be treated in the same manner, but it can be removed without compromising any of the considerations in the paper. 
 
    Capillaries are characterized by a slender geometry. This means that, not only $R\ll \mathrm{diam}(\Omega)$, but also $R\ll L$, where $L$ denotes the total length of a considered capillary network. It is hence common in literature {(see, among others,  \cite{Zunino2019, Koch2020, Koppl, Fritz, Kuchta2021})} to reduce the 3D-3D problem, in which $\D^k$ and $\Sigma^k$ are coupled through some interface condition imposed on $\Gamma^k$, to a 3D-1D coupled problem, where the capillaries are identified with their centerlines and the interstitial volume is extended to the whole tissue sample. This enables a strong reduction of the meshing complexity and of the computational cost. In this spirit we hence assume that, in each time interval $\mathcal{I}_k$, the original network $\Sigma^k$ is covered by a collection of slender and {straight} cylindrical vessels $\Sigma_i^k$, each with lateral surface $\Gamma_i^k$ and centerline $\Lambda^k_i=\lbrace \bm{\lambda}_i^k(s), ~s\in (0, S_i^k)\rbrace$, for $i\in Y^k$, {where $Y^k$ indexes the vessel segments in $\mathcal{I}_k$. Vessel centrelines may connect at endpoints $\lbrace\bm{x}_b\rbrace_{b \in B^k}$, with $B^k$ denoting the set of connection points. The subset $Y_b^k\subset Y^k$ contains the indices of segments joined at $\bm{x}_b$, and $S_{i,b}=\lbrace 0,S_i^k\rbrace$ is the curvelength coordinate such that $\bm{\lambda}_i^k(S_{i,b})=\bm{x}_b$, $i \in Y_b^k$. In the derivation of the 3D-1D problems, the domain $\Sigma^k$ is thus replaced, by $\bigcup_{i\in Y^k} \Sigma_i^k$, and then identified with $\Lambda^k=\bigcup_{i\in Y^k}\Lambda_i^k \cup \lbrace\bm{x}_b\rbrace_{b \in B^k} $; further, {$\Lambda_{in}$, $\Lambda_{out}$ and $\Lambda_d^k$ correspond to the sets of the centers of the sections in $\partial \Sigma_{in}$, $\partial \Sigma_{out}$ and $\partial \Sigma_d^k$}, respectively. We remark that, in general, according to the above definition, $\Sigma^k\neq \bigcup_{i \in Y^k}\Sigma_i^k$, but, as $R$ is small such a difference can be considered negligible. 
    
    Let $q(\xx,t)$ and $\tilde{q}(\xx,t)$ be two generic quantities defined respectively in $\D^k$ and $\Sigma^k$. For all $ t \in \I_k$, $i \in Y^k$ and for any $s \in (0,S_i^k)$ we assume, using cylindrical coordinates aligned with $\Lambda_i^k$, that 
    \begin{equation}
\tr{q}{\Gamma_i^k}(\theta,s,t;R)=\check q_i(s,t)\quad \forall \theta \in [0,2\pi) 
\label{qcheck}
\end{equation}
\begin{equation}
\tr{\tilde q}{\Sigma_i^k}(r,\theta,s,t)=\hat q_i(s,t)\quad \forall r \in[0,R],~\forall \theta \in [0,2\pi).\label{qhat}
\end{equation}
In other words, due to the slender geometry of the vessel segments, we are assuming that the variations of the quantities of interest on their cross sections are negligible. Specifically, in \eqref{qcheck}, the restriction of a quantity defined in the interstitial volume to the lateral surface of a vessel segment is assumed to vary only with the curvilinear coordinate $s$. This allows us to define a one-dimensional counterpart of trace variables without actually defining a 3D-1D trace operator. Similarly, in \eqref{qhat}, quantities defined within the vessel network are considered to have uniform cross-sectional profiles and therefore also depend uniquely on $s$.  We refer to \cite{BGS3D1Ddisc} for a rigorous derivation of a well-posed 3D-1D coupling resorting to the above assumptions. 

Referring to the symbols introduced in \eqref{qcheck}-\eqref{qhat}, in the following we adopt the notation $$\check{q}=\prod_{i \in Y^k}\check q_i \quad \text{ and } \quad \hat{q}=\prod_{i \in Y^k}\hat q_i $$ and, as the capillaries are now identified with their centerlines, we will assume that $\mathcal{D}^k\equiv \Omega$.
     
\section{Mathematical model} \label{sec:model}

	Let us consider a cell aggregate modeled as a mixture of cells, interstitial fluid and extracellular matrix (ECM) components, whose respective volume fractions are denoted by $\phic$, $\phil$ and $\phiecm$. We assume the extracellular matrix to be rigid and inert, hence we disregard the evolution of $\phiecm$. For any $k=0,...,K$ we have
	\begin{align}\phic^k+\phil^k=1-\phiecm=\phimax \label{eq:sum_phi} \quad \text{in } \Omega ,
	\end{align} 
  where $\phic^k=\phic(t_k)$ and $\phil^k=\phil(t_k)$. We here assume that tumor cells can not penetrate the capillary wall and enter blood flow. On the other hand, the fluid phase is present both in the mixture, where it represents the interstitial fluid, and in the vessel network, where it models blood.
	
	 We aim at describing the variation of $\phic$ and $\phil$ in response to the oxygen and fluid supplied by blood flow during angiogenesis. The growth of the capillary network is directly promoted by tumor cells, through the production of a vascular endothelial growth factor (VEGF).
	The process is described by means of four partial differential equation problems (PDEs) and an ordinary differential equation (ODE). The PDEs describe the variation in time and space of the volume fraction of tumor cells, of the fluid pressure of the cell aggregate and blood, of oxygen and of VEGF concentration, whereas the ODE is used to model the growth of the capillaries.
	
	The starting point is an equilibrium configuration for all the involved quantities. Vessel growth is induced by the presence of VEGF: tip cells respond chemotactically to VEGF gradients, moving towards regions in which the concentration of the chemical is higher. The growth of the network in a given time interval breaks the equilibrium, causing the evolution of the quantities of interest over the same time interval. 
	First of all the new distribution of $\phic$ in the interstitial volume is computed by solving a non-linear PDE, depending on oxygen availability. The liquid volume fraction $\phil$ is then obtained from \eqref{eq:sum_phi} and used to compute the pressure of the aggregate in the interstitial volume and blood pressure in the capillaries. Finally, fluid velocities are computed by post processing the pressure variables and used in the transport term in the oxygen and VEGF problem. Here, the magnitude of the VEGF source is directly proportional to the tumor cell volume fraction and is also modulated by oxygen availability, such that capillaries are preferentially attracted to the most hypoxic regions.
    
   In the following the different PDE problems are described in detail: first we derive the model problem describing the cell aggregate (tumor cell, liquid phase and ECM), then we move to the oxygen and VEGF concentration problems, ending with the ODE for vessel growth. 
		\subsection{Cell aggregate balance equations}
		 Under the assumption of rigid and inert extracellular matrix, the cell aggregate dynamics in $\Omega$ is described by the variation in time and space of the volume fractions of tumor cells and liquid phase. In particular, for any $t \in \I_k$, mass and momentum balance for all phases in the mixture, together with the saturation constraint, must be satisfied. Hence,
		\begin{align}
		&\phic+\phil=1-\phiecm=\phimax\label{eq:aggregate1}\\
		&\frac{\partial \phic}{\partial t}+\nabla \cdot(\phic \bm{v}_c)=\Gamma_c\label{eq:aggregate2}\\
		&\frac{\partial \phil}{\partial t}+\nabla \cdot(\phil \bm{v}_l)=\Gamma_l\label{eq:aggregate3}\\
		&\nabla \cdot\mathbb{T}_c-\phic\nabla p+\bm{m}_{cl}+\bm{m}_{cm}=0\label{eq:aggregate5}\\
		&-\phil\nabla p+\bm{m}_{lc}+\bm{m}_{lm}=0.\label{eq:aggregate6}
		\end{align}
		Here $\bm{v}_c$ and $\bm{v}_l$ are respectively the tumor cell and the liquid phase velocities, $p$ is introduced as a Lagrange multiplier due to the saturation constraint \eqref{eq:aggregate1} and is
classically identified with the interstitial pressure \cite{Preziosi, GiversoGrillo}, while 
		\begin{align}
		&\bm{m}_{cl}=\phil\phic\mu\mathbb{K}(\phil)^{-1}(\bm{v}_l-\bm{v}_c)=-\bm{m}_{lc},\label{eq:mcl}\\
		&\bm{m}_{cm}=-\mathbb{M}(\phic,\phil)^{-1}\bm{v}_c,\\
		&\bm{m}_{lm}=-\phil\phiecm\mu\mathbb{K}(\phil)^{-1}\bm{v}_l\label{eq:mlm}
		\end{align}
	 describe the interaction forces between the constituents of the mixture, with $\mathbb{M}$ and $\mathbb{K}$ being the motility and permeability tensors respectively, and $\mu$ denoting fluid viscosity.
	 
	 We assume the tumor cell phase to behave as an elastic fluid, characterized by an isotropic stress response governed by a nonlinear pressure-like function $\Sigma(\phi_c)$ depending on the cell volume fraction. Accordingly, the Cauchy stress tensor can be written as
	 \begin{equation}
	 \mathbb{T}_c(\phic)=-\Sigma(\phic)\mathbb{I}=-E\frac{\phic-\phiz}{\phimax-\phic}\phic\mathbb{I},\label{eq:cauchy}
	 \end{equation}
	 where $E$ is the Young modulus, $\phiz$ is the reference volume fraction at which the cellular aggregate experiences zero stress, and $\mathbb{I}$ denotes the identity tensor. 
Denoting by $\Delta p$ the pressure jump between the arterial and the lymphatic compartments, and under the biologically consistent assumption that $\frac{\Delta p}{E} \ll 1$, it can be shown \cite{GiversoGrillo} that the pressure-drop contribution in Eq.~\eqref{eq:aggregate5} can be neglected and that the interaction between the cells and the liquid phase is negligible compared to the interaction between the cells and the extracellular matrix, i.e., $|\bm{m}_{cl}| \ll |\bm{m}_{cm}|$.
Furthermore, under the biologically admissible assumption that $|\bm{v}_l| \gg |\bm{v}_c|$, the cell-velocity contribution in Eq.~\eqref{eq:aggregate6} can be disregarded \cite{GiversoGrillo}.
Finally, we assume both the permeability and the motility tensors to be constant and isotropic, namely $\mathbb{K} = \kappa \mathbb{I}$ and $\mathbb{M} = M \mathbb{I}$. Consequently, Eqs.~\eqref{eq:aggregate1}--\eqref{eq:aggregate6} can be rewritten as
 \begin{align}
  &\frac{\partial \phic}{\partial t}-\nabla \cdot(M \phic \Sigma'(\phic)\nabla\phic)=\Gamma_c\label{eq:tumcell}\\
 &\phic+\phil=\phimax\label{eq:algebraic}\\
 &\frac{\partial \phil}{\partial t}-\nabla \cdot \Big(\frac{\phil}{1-\phil}\frac{\kappa}{\mu}\nabla p\Big)=\Gamma_l.\label{eq:pressure}
 \end{align}
We now provide further details on the tumor cell equation \eqref{eq:tumcell} and on the pressure equation \eqref{eq:pressure}, concerning source terms and boundary and initial conditions. Equation \eqref{eq:pressure} is referred to as pressure problem as $\phil$ is directly obtained by the algebraic equation \eqref{eq:algebraic} before solving \eqref{eq:pressure}.
 \paragraph{The tumor cells problem}
      Cell proliferation is influenced by multiple chemical and mechanical factors, including the amount of cells already present in the tissue, the availability of free space, and the supply of nutrients such as glucose and oxygen. In the present model, we assume that oxygen concentration, denoted by $c$, acts as the limiting factor. Accordingly, we adopt the following expression for the reaction term
\begin{equation}
\Gamma_c(\phi_c,c) = \phi_c S_c(\phi_c,c),
\end{equation}
with
\begin{equation}
S_c(\phi_c,c) =
{\frac{\gamma}{c_{\mathrm{ref}}}((\phi_{\max} - \phi_c)c
-c_{\mathrm{ref}})_+}
\label{eq:Sc}
\end{equation}
{where $\gamma$ denotes the intrinsic tumor cell doubling rate, $(\cdot)_+$ represents the positive part operator, and $c_{\mathrm{ref}}$ is the threshold {volumetric} oxygen concentration required for cell proliferation.}
{In the present work, tumor cell penetration through the capillary wall is not considered; accordingly, no evolution equation is introduced for tumor cells inside the capillaries.} 
{Considering the Cauchy stress tensor given by Equation \eqref{eq:cauchy}, we can define}
		\begin{equation}
		F_c(\phi_c)=M\phic\Sigma'(\phic)=EM(\phic)\phic\frac{\phimax(2\phic-\phiz)-\phic^2}{(\phimax-\phic)^2},
		\end{equation}
	Then, for $t \in \I_k$ we look for $\phic$ as the solution of
			\begin{align}
		&\frac{\partial \phic(\xx,t)}{\partial t}-\nabla \cdot (F_c(\phic(\xx,t))\nabla \phic(\xx,t))-S_c(\phic(\xx,t),c(\xx,t))\phic(\xx,t)=0 &\xx \in\Omega\label{eq:phic1}\\
		&F_c(\phic(\xx,t))\nabla \phic(\xx,t)\cdot \bm{n}(\xx)=0 &\xx \in\partial\Omega,\label{eq:phic2}
		\end{align}
		where $\bm{n}$ denotes the unit normal vector to $\partial \Omega$. {The zero-flux boundary condition \eqref{eq:phic2} expresses the fact that tumor cells are not allowed to enter or leave the domain.} 
		
		Concerning the initial conditions, for the moment we just assume the values of $\phic$ and $c$ at time $t_{k-1}$ to be given data (see Section \ref{subsec:discrete_tum_cell} for details). 
        
		\paragraph{The pressure problem}
		{The pressure problem governs fluid motion within the mixture and must be coupled with the flow in the vessel network.} The source term in the liquid-phase mass balance \eqref{eq:pressure} hence accounts for {liquid absorption by proliferating tumor cells} and for the fluid flux from the embedded vascular and lymphatic networks}. Therefore, it takes the following form 
		$$\Gamma_l(\phic,\phil,p,\hat p)=-\phic S_c(\phic,c)-\Gamma_{LS}(p)+\Gamma_\Lambda(p,\hat p),$$ 
        The last two terms account for the liquid exchanges with the lymphatic system and the capillary network respectively, with $\hat p$ denoting blood pressure along the one dimensional capillary network $\Lambda^k$. In particular, 
        {we model the lymphatic system as a homogenized source/sink of fluid within the network, activated when the interstitial pressure is below or above the lymphatic pressure $p_{LS}$, assumed constant. Accordingly, we set}
		\begin{equation*}
		\Gamma_{LS}(p)=\phil\beta_p^{LS}\frac{S}{V}(p-p_{LS}),
		\end{equation*}
         being $\beta_p^{LS}$ the permeability of the lymphatic wall, and $\frac{S}{V}$ the surface area of the lymphatic vessels per unit tissue volume.
{As for the vascular network, we retain an explicit description of the individual vessels. Accordingly, we introduce a source/sink term concentrated on the vessels, of the form}
$$\Gamma_\Lambda(p,\hat{p})=\sum_{i \in Y^k}\hat{f}_p^i(s,t)\delta_{\Lambda_i^k},$$
where {$\delta_{\Lambda_i^k}$ denotes the Dirac delta distribution centered on $\Lambda_i^k$.}
{The fluid flux $\hat{f}_p^i$ exchanged at the vessel wall is derived from a modified Starling law, which distinguishes filtration and absorption regimes, to accounts for the liquid volume fraction within the tissue and for the fact that the vessels are fully filled with fluid. Indeed, the liquid volume fraction in the capillary network is assumed to be equal to one. Accordingly, we set}
\begin{equation}
	\hat{f}_p^i:=\begin{cases}
	2\pi R\beta _p(\hat p_i-\check p_i-\Delta p_{onc}) & \text{ if }\hat p_i>\check p_i+\Delta p_{onc}\\2\pi R\beta _p~\check{\phi}_{l,i}~(\hat p_i-\check p_i-\Delta p_{onc}) & \text{ if }\hat p_i<\check p_i+\Delta p_{onc}
	\end{cases}\quad s \in[0,S_i^k]. \label{eq:fpi}
	\end{equation}
For the $\check{(\cdot)}$ and $\hat{(\cdot)}$ notation we refer to {Equations \eqref{qcheck}-\eqref{qhat}}.
{The Starling law in \eqref{eq:fpi} assumes that the fluid flux across the capillary wall is proportional to the pressure jump through a positive scalar coefficient $\beta_p$, representing the permeability of the capillary wall. The pressure jump accounts for the difference between the fluid pressures in the tissue and in the vessel network—both model unknowns—as well as for the oncotic pressure contribution, which is treated as constant \cite{cattzun}. 
Indeed, the oncotic pressure difference $\Delta p_{onc}$ is determined by the solute concentration difference across the vessel wall. Since the dominant contribution to oncotic pressure is due to plasma proteins, primarily albumin, whose concentration is assumed constant, the term $\Delta p_{onc}$ is also taken to be constant.}

		 Then, given $\phic$ from \eqref{eq:phic1}-\eqref{eq:phic2} and $\phil$ from \eqref{eq:algebraic}, we look for the pressure $p$ in $\Omega$ for $t \in \I_k$ as the solution of
		\begin{align}
&\frac{\partial \phil(\xx,t)}{\partial t}-\nabla \cdot \Big(\frac{\phil(\xx,t)}{1-\phil(\xx,t)}\frac{\kappa}{\mu}\nabla p(\xx,t)\Big)+\phil(\xx,t)\beta_p^{LS}\frac{S}{V}(p(\xx,t)-p_{LS})=\nonumber\\&\hspace{3cm}-\phic(\xx,t) S_c(\phic(\xx,t),c(\xx,t))+\sum_{i \in Y^k}\hat f_i^p\delta_{\Lambda_i^k}\label{eq:pr3D} \hspace{1.1cm} \xx \in \Omega\\
&\Big(\frac{\phil(\xx,t)}{1-\phil(\xx,t)}\frac{\kappa}{\mu}\nabla p(\xx,t)\Big)\cdot \bm{n}(\xx)=0 \hspace{2.5cm}\xx \in \partial \Omega\label{eq:p2}.
		\end{align}
		The system is closed by introducing a PDE problem describing the blood pressure within each vessel segment. {Specifically, blood flow inside the capillaries is modeled as an incompressible viscous fluid governed by Poiseuille law for steady laminar flow \cite{Baxter_1, cattzun0}, coupled with the exchange flux \eqref{eq:fpi} at the vessel wall. In the 3D-1D framework, this yields the following problem \cite{BGGPS2023},} for $t \in \I_k$:
		\begin{align}
		-&\frac{\pi R^4}{8 \mu}\frac{d^2\hat p_i(s,t)}{ds^2}=-\hat{f}_p^i(s,t)& s \in (0,S_i^k),\label{eq:pr1D}\\
		&\hat{p}=\hat{p}_{in} &\text{ in }\Lambda_{in}\\
		&\hat{p}=\hat{p}_{out} &\text{ in }\Lambda_{out}\\
		&\frac{d\hat{p}}{ds}=0 &\text{ in }\Lambda_{d}^k\\
		&\sum_{j \in Y_b} \frac{\partial \hat{p}_j}{\partial s}(S_{j,b},t)=0 & \forall b\in B^k\\
		&\hat{p}_i(S_{i,b},t)=\hat{p}_j(S_{j,b},t) &\forall i\neq j \in Y_b,~\forall b \in B^k,\label{pr1D_end}
		\end{align}	
		where the last two equations express respectively flux balance and pressure continuity at bifurcation points. We remark that the flow through the lateral surface of each vessel, $\hat{f}_p^i$ behaves as a distributed source term in the 1D equation~\eqref{eq:pr1D}, balanced by the concentrated source term in the 3D equation~\eqref{eq:pr3D}.
		
		Once the fluid pressure is obtained in $\Omega$ and for all $\Lambda_i^k$, the fluid velocity can be computed as
		\begin{align}
&\bm{v}_l(\xx,t)=-\frac{1}{1-\phil(\xx,t)}\frac{\kappa}{\mu}\nabla p(\xx,t), &\quad \bm{x} \in \Omega,~t \in \I_k\\
&\hat{v}_{l,i}(s,t)={-}\frac{R^2}{8\mu}\frac{d\hat p_i(s,t)}{ds},&\quad s\in(0,S_i^k),~t \in \I_k
		\end{align}
		
\subsection{The oxygen concentration {model}}
		Let $t\in \I_k$ and let us denote by $c(\bm{x},t)$ the oxygen concentration {in the liquid phase within $\Omega$} and by $\hat{c}_i(s,t)$ the oxygen concentration in each $\Lambda_i^k$, $i \in Y^k$. Let us also define {the flux exchange at the vessel wall}
		\begin{equation}
		\hat{f}_c^i(s,t)=\begin{cases}
		2\pi R\beta _c(\hat c_i(s,t)-\check c_i(s,t))& \text{ if }\hat c_i>\check c_i\\2\pi R\beta _c~\check{\phi}_{l,i}~(\hat c_i(s,t)-\check c_i(s,t)) & \text{ if }\hat c_i<\check c_i
		\end{cases}
		\end{equation}
 where $\beta_c$ is the permeability to oxygen of the vessel wall and $\check{c}_i=\tr{c}{\Gamma_i^k}$. {Also in this case, we distinguish between oxygen absorption and release at the vessel wall to account for the fact that nutrient uptake from the tissue is proportional to the local liquid volume fraction.} The 3D-1D reduced problem describing the diffusion and transport of oxygen in the interstitial volume and in the capillary network for $t \in \I_k$ is basically the same that was proposed in \cite{BGGPS2023}, with the difference that the oxygen is here dissolved in the liquid phase:
	\begin{align}
	&\phi_l(\xx,t) \frac{\partial c(\xx,t)}{\partial t}-\nabla \cdot(\phil(\xx,t) D_c\nabla c(\xx,t))+\bm{v}_l(\xx,t)\phil(\xx,t)\nabla c(\xx,t)\nonumber\\[-0.5em]&\hspace{4cm}+ m_c \phil(\xx,t)c(\xx,t)(\phimax-\phil(\xx,t))=\sum_{i \in Y^k}\hat f_i^c\delta_{\Lambda_i^k}&\xx \in \Omega\\
	&\phil(\xx,t)D_c\nabla c(\xx,t)\cdot \bm{n}(\bm{x})=0 &\xx \in \partial \Omega\\[1em]
	&\pi R^2\frac{\partial \hat c(s,t)}{\partial t}-\pi R^2 \tilde{D}_c\frac{\partial^2\hat c(s,t) }{\partial s^2}+\hat{v}_l(s,t)\frac{\partial \hat c(s,t)}{\partial s}=-\hat{f}_c^i(s,t)&s \in (0,S_i^k)\\
		&\hat{c}=\hat{c}_{in} &\text{ in }\Lambda_{in}\\
	&\frac{\partial\hat{c}}{\partial s}=0 &\text{ in }\Lambda_{out}\cup\Lambda_{d}^k\\
	&\sum_{j \in Y_b} \frac{\partial \hat{c}_j}{\partial s}(S_{j,b},t)=0 & \forall b\in B^k\\
	&\hat{c}_i(S_{i,b},t)=\hat{c}_j(S_{j,b},t) &\hspace{-1cm}\forall i\neq j \in Y_b,~\forall b \in B^k
	\end{align}
	Parameters $D_c$ and $\tilde{D}_c$ are positive scalars denoting oxygen diffusivity respectively
	in {the tissue and in the capillary network,} whereas $m_c$  is the rate at which oxygen is metabolized by tumor cells. For what concerns the initial condition, for $k=1$ we set $c(\bm{x},t_{0})=c^0(\bm{x})$ in $\Omega$ and $\hat c(\bm{x},t_{k-1})=\hat c^0(\bm{x})$ in $\Lambda^0$. For $k>1$ oxygen concentration at time $t_{k-1}$ is directly available in $\Lambda^{k-1}$ from the final concentration computed in the time interval $\mathcal{I}_{k-1}$, while we choose to set $\hat{c}(\xx,t_{k-1})=0$ for $\xx \in \Lambda^k\setminus \Lambda^{k-1}$, {i.e. for the capillary segments generated during $\I_k$.}
    
	\subsection{The chemotactic growth factor {model}}
\begin{figure}
		\centering
		\includegraphics[width=0.5\linewidth]{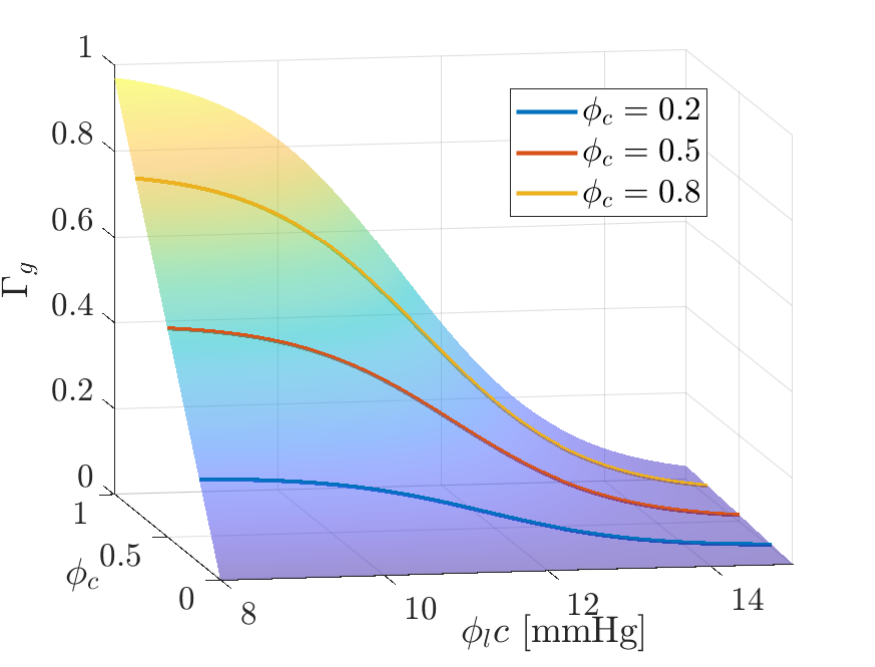}
		\caption{VEGF source term $\Gamma_g(\phic,c)$ for $G=1$,  $c^\star=11.5~ \mathrm{mmHg}$ and $b=11.5$.}
		\label{fig:fvegf}
	\end{figure}

	The third PDE problem models the distribution of VEGF in the tissue. The consumption of VEGF by endothelial cells is intended in terms of receptor mediated binding at the vessel wall, meaning that there is no VEGF flux through the capillary wall. For this reason no equation is defined for the concentration of VEGF inside the capillaries.
   {Let $g$ denote the VEGF concentration in the liquid phase within $\Omega$, and let $\check{g}_i = \tr{g}{\Gamma_i^k}$. The absorption of VEGF at the vessel wall can then be represented in the following form}
	\begin{equation*}
	\hat{f}_i^g(s,t)=2\pi R\tilde{\sigma}\check{\phi}_{l,i}(s,t)\check{g}_i(s,t) \,
	\end{equation*} 
    {being $\tilde{\sigma}$ the rate of VEGF consumption by endothelial cells.}
    {The evolution of the VEGF concentration in $\Omega$ for $t \in \I_k$, accounting for chemical diffusion, fluid advection, and production/absorption mechanisms, is governed by the following advection-diffusion-reaction equation:}
	\begin{align}
	&\phil(\xx,t)\frac{\partial g(\bm{x},t)}{\partial t}-\nabla \cdot \big(\phil(\xx,t)D_g \nabla {g}(\bm{x},t)\big)+\bm{v}_l(\bm{x},t)\phil(\xx,t)\cdot\nabla g(\bm{x},t)+ \nonumber\\[-0.4em]
&\hspace{3.5cm}+\sigma \phil(\xx,t)g(\bm{x},t)=\Gamma_g(\phi_c(\bm{x},t),c(\bm{x},t))-\sum_{i \in Y^k}\hat{f}_i^g(\check g)\delta_{\Lambda_i^k}, \quad \bm{x} \in \Omega
\label{eq_g_strong}
\end{align}
Parameters $D_g$ and  $\sigma$ are positive scalars denoting, respectively, the diffusivity of the VEGF and its natural decay rate. {Furthermore, compared with \cite{BGGPS2023}, the VEGF production term $\Gamma_g(\phi_c(\bm{x},t),c(\bm{x},t))$ explicitly depends on both the tumour cell volume fraction and the local oxygen concentration. Consistently with the underlying biological mechanisms, VEGF production by cancer cells is mainly triggered in hypoxic regions, thereby promoting capillary growth toward these areas and improving vascular perfusion.} For this reason we define a sigmoidal forcing term 
	$$\Gamma_g(\phic,c)=G\phic\Big(1-\frac{1}{1+e^{b(1-{\phil}c/c^\star)}}\Big)$$
	where $G$ is a positive scalar denoting the VEGF production rate by tumor cells, $c^\star$ is a reference oxygen concentration {per unit of tissue volume} for which $\Gamma_g=\frac{1}{2}G\phi_c$ and $b$ is a dimensionless constant allowing us to control the steepness of the sigmoid in the surroundings of the inflection point. The behavior of function $\Gamma_g$ is reported in Figure $\ref{fig:fvegf}$ for $G=1$ and $c^\star=11.5~ \mathrm{mmHg}$. The shape of $\Gamma_g$ is chosen to promote VEGF production predominantly in hypoxic tumor regions, 

    Equation~\eqref{eq_g_strong} is complemented by a zero-flux boundary condition
\begin{equation}
\phil(\xx,t) D_g \nabla g(\xx,t)\cdot \bm{n}(\xx) = 0
\qquad \xx \in \partial \Omega \, .
\end{equation}
In analogy with the other model variables, the initial condition for $k=1$ coincides with the equilibrium configuration, namely $g(\xx,t)=g_0(\xx)$, whereas for $k>1$ it is inherited from the final VEGF concentration computed at the previous time step.

	\subsection{The capillary growth {model}}
	The growth of the capillary network is modeled as in \cite{BGGPS2023}. For this reason we here recall only the main aspects, referring the reader to \cite{BGGPS2023} for further details.
	
	The growth of the capillary network is modeled by a tip-tracking approach, in which sprout extension is rendered as a displacement of the tip cells. The position $\bm{x}_{P}$ of a generic tip cell evolves according to 
	\begin{equation}
	\frac{d\bm{x}_P}{dt}=\bm{w}(g(\bm{x}_P,{t_{k}}),\bm{x}_P)\label{grow}
	\end{equation}
	with $\bm{w}$ denoting the tip velocity and defined, according to \cite{sun}, as
	\begin{equation}
	\bm{w}(g,\bm{x})=\begin{cases}\cfrac{l_e}{t_c(g)}\cfrac{\bm{K}_{ECM}(\bm{x})\nabla g}{||\bm{K}_{ECM}(\bm{x})\nabla g||} &\text{ if }g\geq g_{lim}\\
	0 &\text{otherwise}.\end{cases} \label{growth_vel}
	\end{equation}
	Parameter $g_{lim}$ represents the minimum VEGF concentration for endothelial cell proliferation, $l_e$ is the endothelial cell length and $t_c$ is a cell cycle division time, modeled as \cite{sun}
	\begin{equation}
	t_c(g)=\tau \left(1+e^{\left(\frac{\bar{g}}{g}-1\right)}\right),\label{tc}
	\end{equation}
	where $\tau$ is a cell proliferation parameter, while $\overline{g}$ is the {VEGF concentration in the liquid phase} at which $t_c=2\tau$. The local orientation of the extracellular matrix (ECM) fibers is modeled through matrix $\bm{K}_{ECM}$, which is assumed constant in time but variable in space. When no data on the ECM orientation is available, the matrix is defined as a random perturbation of the identity matrix (see \cite{BGGPS2023}).
	
\begin{figure}
	    \centering
	    \includegraphics[width=0.5\linewidth]{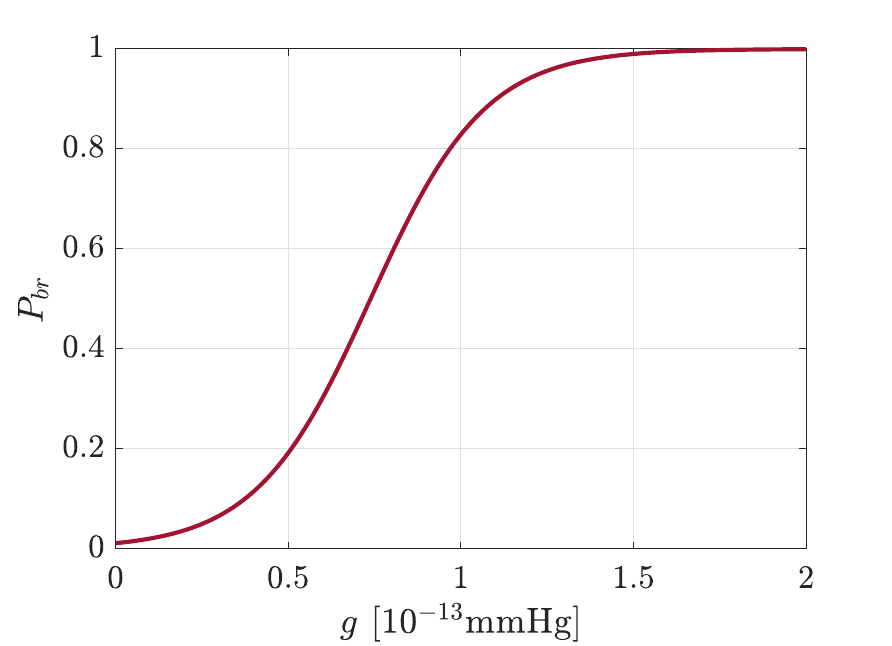}
	    \caption{{Example of branching probability function used in the simulations. In this case $P_{br} (g)= \frac{1}{1 + e^{-a(1 - d)}}$ with $a$ and $d$ such that $P_{br}(g_{br})=0.99$ and $P_{br}(g_{lim})=0.05$, for $g_{br}=1.5\cdot 10^{-13}$ mmHg  and $g_{lim}=0.25 \cdot 10^{-13}$ mmHg. }}
	    \label{fig:Pbr}
	\end{figure}

	Branching can occur when the age of the sprout is greater than a threshold age $\tau_{br}$, and the component of $\bm{w}$ perpendicular to the parent sprout orientation is big enough ($||\bm{w}_\perp||\geq \alpha_{br}^w||\bm{w}||$). To ensure that branching occurs predominantly at high VEGF concentrations, we introduce a sigmoidal branching probability function 
	$P_{br}(g)$. The analytical form of this function is chosen on a qualitative basis, in order to ensure very little branching {in case the VEGF concentration in the liquid phase is} close to $g_{lim}$ and to have control over the VEGF concentration $g_{br}$ corresponding to a branching probability close to 100\%. An example of $P_{br}(g)$ is reported in Figure \ref{fig:Pbr}, {but different choices are of course possible}.

	If branching occurs, the directions of the two resulting sprouts are are obtained as
	\begin{equation}
	\bm{w}_{1,2}=\bm{w}\pm\frac{\bm{w}_{\Pi}}{||{\bm{w}_{\Pi}}||}d_{br}\label{eq:branching_direction}
	\end{equation} with $d_{br}$ denoting the distance by which the new sprout tips will be separated.
	
	\subsection{Summary of model equations}	
	We report here a concise summary of the model equations. {For the sake of compactness, and to better highlight the possible non-linear relations among the unknowns, we suppress in the notation the dependence on the spatial coordinate $\xx$ and on the time $t$.}
	\begin{align}
	&\begin{cases}
	&\frac{\partial \phic}{\partial t}-\nabla \cdot (F_c(\phic)\nabla \phic)-S_c(\phic,c)\phic=0 \\
	&\phic+\phil=\phimax
	\end{cases}\hspace{5.2cm} \text{in } \Omega\times \I_k\label{eq:modelsummary_1}\\[0.5em]
	&\begin{cases}
\frac{\partial \phil}{\partial t}-\nabla \cdot \Big(\frac{\phil}{1-\phil}\frac{\kappa}{\mu}\nabla p\Big)+\phil\beta_p^{LS}\frac{S}{V}(p-p_{LS})+\phic S_c(\phic,c)=\sum_{i \in Y^k}\limits\hat f_i^p(\hat{p},\check{p})\delta_{\Lambda_i^k} \quad \text{in } \Omega\times \I_k\\
-\frac{\pi R^4}{8 \mu}\frac{d^2\hat p_i}{ds^2}=-\hat{f}_p^i(\hat{p},\check{p}) \hspace{6.3cm} 	\text{ in }(0,S_i^k)\times \I_k ~\forall i \in Y^k
	\end{cases}\label{eq:modelsummary_2}\\[0.5em]
		&\begin{cases}\phi_l \frac{\partial c}{\partial t}-\nabla \cdot(\phil D_c\nabla c)+\bm{v}_l(p)\phil\cdot\nabla c+ m_c \phil c(\phimax-\phil)=\sum_{i \in Y^k}\limits\hat f_i^c(\hat{c},\check{c})\delta_{\Lambda_i^k}\quad\text{ in } \Omega\times \I_k\\
			\pi R^2\frac{\partial \hat c}{\partial t}-\pi R^2 \tilde{D}_c\frac{\partial^2\hat c }{\partial s^2}+\hat{v}_l(\hat{p})\frac{\partial \hat c}{\partial s}=-\hat{f}_c^i(\hat{c},\check{c}) \hspace{3.3cm}	\text{ in }(0,S_i^k)\times \I_k ~\forall i \in Y^k\end{cases}\label{eq:modelsummary_3}\\[0.5em]&
				\phil\frac{\partial g}{\partial t}-\nabla \cdot \big(\phil D_g \nabla {g}\big)+\bm{v}_l(p)\phil\cdot\nabla g+\sigma \phil g=\Gamma_g(\phi_c,c)-\sum_{i \in Y^k}\hat f_i^g(\check g) \delta_{\Lambda_i^k}\hspace{1.1cm} \text{ in } \Omega\times \I_k
				\label{eq:modelsummary_4}
	\end{align}
	The system is closed with the following boundary conditions:
	\begin{align*}
	&F_c(\phic)\nabla \phic\cdot \bm{n}=0, \quad \Big(\frac{\phil}{1-\phil}\frac{\kappa}{\mu}\nabla p\Big)\cdot \bm{n}=0 &\text{ in } \partial\Omega\times \I_k\\&\phil D_c\nabla c\cdot \bm{n}=0 , \quad\phil D_g\nabla g\cdot \bm{n}=0  &\text{ in } \partial\Omega\times \I_k\\
	&\hat{p}=\hat{p}_{in},\quad \hat{c}=\hat{c}_{in}  &\text{ in }\Lambda_{in}\times \I_k\\
	&\hat{p}=\hat{p}_{out}, \quad \hat{c}=0 &\text{ in }\Lambda_{out}\times \I_k\\
	&\frac{d\hat{p}}{ds}=0, \quad \frac{d\hat{c}}{ds}=0 &\text{ in }\Lambda_{d}^k\times \I_k
	\end{align*}
	and the conditions at capillary network connection points:
	\begin{align*}
	&\sum_{j \in Y_b} \frac{\partial \hat{p}_j}{\partial s}(S_{j,b},t)=0, \quad \sum_{j \in Y_b} \frac{\partial \hat{c}_j}{\partial s}(S_{j,b},t)=0 & \forall b\in B^k\\
	&\hat{p}_i(S_{i,b},t)=\hat{p}_j(S_{j,b},t), \quad \hat{c}_i(S_{i,b},t)=\hat{c}_j(S_{j,b},t) &\hspace{-1cm}\forall i\neq j \in Y_b,~\forall b \in B^k
	\end{align*}
	For what concerns initial conditions, we set:
	\begin{align*}
	\phic=\phic^0 , \quad p=p^0, \quad c=c^0, \quad &g=g^0&\text{ in } \Omega \times \lbrace t=0\rbrace\\
	&\hat{c}=0&\text{ in }\Lambda^k\setminus \Lambda^{k-1} 
	\end{align*}
    
	\section{Problem discretization and solving strategy} \label{sec:numerics}
	The growth of the capillary network and the evolution of the quantities of interest are decoupled by means of a continuous-discrete hybrid tip-tracking strategy. Under this approach, Equation \eqref{grow} is solved using a forward Euler scheme in the time interval $\I_k$. Once the capillary geometry has been updated, we consider the evolution of the quantities of interest for $t\in \I_k$ on the fixed geometry. In particular, we discretize \eqref{eq:modelsummary_1}-\eqref{eq:modelsummary_4} using linear finite elements in space and a backward Euler scheme in time. Non-linearities are handled by the Newton scheme, whereas linear 3D-1D coupled problems are recast into PDE constrained optimization problems, following the method introduced in \cite{BGS3D1Ddisc}, and previously used for angiogenesis simulations in \cite{BGGPS2023}.  We here focus mainly on the novelties introduced with respect to \cite{BGGPS2023}, hence on the tumor cell and liquid phase volume fraction equations. More details on the optimization based approach can be found in the above mentioned references.
	
		Let us consider a tetrahedral mesh $\mathcal{T}_h$ of $\Omega$, completely independent from the position of the vessel network and fixed in time. We denote by $N$ the number of degrees of freedom, which is the same for all unknowns defined in $\Omega$ and equal to the number of nodes in the grid, as Neumann boundary conditions are imposed. Given a generic quantity $q=q(\xx,t)$, that could correspond to $\phic$, $p$, $c$, or $g$, we define its discrete approximation as 
	$$Q(\xx,t)=\sum_{j=1}^NQ_j(t)\eta_j(\xx),$$ where $\lbrace \eta_j\rbrace_{j=1}^N$ are linear Lagrangian basis functions. Let us observe that the considered basis functions are continuous on $\Omega$. We can hence compute trace variables as $\check Q(\xx,t)=\sum_{j=1}^N\check Q_j(t)\tr{\eta_j}{\Lambda}(s).$
	 
	Variables that are defined directly on the 1D capillary network (namely $\hat p$ and $\hat c$) are semi-discretized in space using piecewise linear finite elements on uniform subpartitions of the straight vessel segments. The conditions at capillary connection points are strongly imposed by defining a single degree of freedom at intersections. Denoting by $\lbrace\hat\eta_j\rbrace_{j=1}^{\hat{N}}$ the set of linear Lagrangian basis functions on $\Lambda^k$, {obtained by assembling the segment-wise basis functions}, we define the approximation of a generic 1D variable $\hat{w}$ as
	$$\hat{W}(s,t)=\sum_{j=1}^{N}\hat{W}_j(t)\hat\eta_j(s,t).$$
	We will have $\hat{N}=\hat{N}_p$ for the 1D fluid pressure and $\hat{N}=\hat{N}_c$ for the 1D oxygen concentration.
		
		In the following, upper case symbols are used to denote either the semi-discrete approximation of a quantity or the column vector collecting the degrees of freedom of its semi-discrete approximation $Q(t)=[Q_1(t),....,Q_N(t)]^T$, $\hat{W}(t)=[\hat{W}_1(t),....,\hat{W}_{\hat{N}}(t)]^T$, the meaning being clear from the context. 
        
        {The remainder of this section is organized as follows. In Section \ref{sec:3D1D} we provide an overview on how general 3D-1D coupled problems are treated within the adopted optimization-based framework. In Section \ref{sec:capillary_update} we give some details on the implementation of the capillary network update. Finally, in Sections from \ref{subsec:discrete_tum_cell} to \ref{subsec:discrete_VEGF} we detail  on the discretization of the PDE models of interest, with particular attention on the time semi-discrete equations and on the treatment of non-linearities. }
	
    \subsection{{Treatment of 3D-1D coupled problems}}
    \label{sec:3D1D}
    The coupled 3D-1D problems for pressure distribution and oxygen concentration in equations \eqref{eq:modelsummary_2} and \eqref{eq:modelsummary_3} are tackled by means of a PDE-constrained domain decomposition strategy. We briefly sketch here the approach, referring to references \cite{BGGPS2023, BGS3D1Ddisc} for all the details. 
    Let us consider, at a given time frame, two coupled linear systems:
    \begin{eqnarray*}
    \begin{cases}
        \bm{A}_\Omega Q + \bm{B}_{\Omega,\Lambda}\left(\check{Q}-\hat{Q}\right)=b_\Omega\\
        \bm{A}_\Lambda \hat{Q} + \bm{B}_{\Lambda,\Lambda}\left(\hat{Q}-\check{Q}\right)=b_\Lambda\\
    \end{cases}
    \end{eqnarray*}
for certain unknown quantities $Q$ defined on a discretization of $\Omega$ and $\hat{Q}$ on a discretization of $\Lambda$. Matrices $\bm{A}_\Omega$, $\bm{A}_\Lambda$ derive from the discretization of the operators set uniquely in the 3D and 1D domains, respectively, whereas matrices $\bm{B}_{\Omega,\Lambda}$ and $\bm{B}_{\Lambda,\Lambda}$ represent the discrete operators that couple the problems on the two domains, being finally $b_\Omega$, $b_\Lambda$ known terms, accounting for forcing terms, boundary conditions and known terms from the previous time-step.

Additional interface variables are then added to decouple the problems, and a cost functional is introduced to express the mismatch between the original and the added variables. The discrete solution is obtained by solving the following optimization problem:
\begin{equation}
\begin{gathered}
    \min J:=\frac{1}{2}\Big(\|\check{Q}-\Psi_\Omega\|_{L^2(\Lambda)}^2+\|\hat{Q}-\Psi_\Lambda\|_{L^2(\Lambda)}^2\Big)\\
    \text{such that:}\\
    \bm{A}_\Omega Q + \bm{B}_{\Omega,\Lambda}\left(\check{Q}-\Psi_\Lambda\right)=b_\Omega,\\
    \bm{A}_\Lambda \hat{Q} + \bm{B}_{\Lambda,\Lambda}\left(\hat{Q}-\Psi_\Omega\right)=b_\Lambda.
\end{gathered}
\label{eq:OPTDISCR}
\end{equation}
The above formulation gives great flexibility in meshing, allowing to introduce independent discretizations for all the involved variables.
The derivation of the optimization formulation for the problems at hand is similar to the one proposed in the mentioned references, where further implementation details and a thorough discussion on the characteristics of the optimization-based approach are also available.

	\subsection{Capillary network update}\label{sec:capillary_update}
We here consider the discrete modeling of the evolution of the capillary network during the time interval $\I_k$. We assume that the initial VEGF concentration $g^{k-1}$ is available either from the previous time step ($k>1$) or from the initial condition ($k=1$). {Let $\Lambda_d^{k-1}$ denote the set of capillary tips at time $t_{k-1}$.} {For each $p \in \Lambda_d^{k-1}$ we introduce a  \textit{virtual} position  $\tilde{\bm{x}}_p^\star$, used to accumulate sub-threshold displacements across time-steps and reset whenever the total displacement exceeds a prescribed threshold $\ell$. Given $\tilde{\bm{x}}_p^\star$ from the previous time-step, we define the \textit{predicted position} as}
\begin{equation}
\tilde{\xx}_p^k=\tilde{\bm{x}}_p^\star+\xx_p^{k-1}+\Delta\I_k \bm{w}(g^{k-1},\xx_p^{k-1})\label{eq:growth_discrete}
\end{equation}
where $\xx_p^{k-1}$ is the position of the $p$-th capillary tip at time $t_{k-1}$ {and $\Delta\I_k \bm{w}$ is the actual {tip elongation during the interval of length $\Delta\I_k= t_k - t_{k-1} $}}. {We then compute the quantity $\Delta\bm{x}_p^k:=\|\tilde{\xx}_p^k-\xx_p^{k-1}\|$,} and define the set of effectively advancing tips as $$\Lambda_{d,\mathrm{eff}}^k=\lbrace p \in \Lambda_d^{k-1}:~\Delta \bm{x}_p^k\geq \ell\rbrace$$ {for a prescribed fixed quantity $\ell>0$, chosen equal to the endothelial cell length $l_c$ {for biological consistency}. For each capillary tip in $ \Lambda_{d,\text{eff}}^k$, we set $\xx_p^k=\tilde{\xx}_p^{k}$, and $\tilde{\bm{x}}_p^\star=0$.
For {$p \in \Lambda_d^{k-1}\setminus \Lambda_{d,\mathrm{eff}}^k$, instead,} the tip position {does not change}, but only the cumulative virtual {position} is {updated}: 
$$\xx_p^k=\xx_p^{k-1}, \quad \tilde{\bm{x}}_p^\star=\tilde{\bm{x}}_p^\star+\Delta\I_k \bm{w}(g^{k-1},\xx_p^{k-1}), \quad \forall p \in \Lambda_{d}^k \setminus \Lambda_{d,\text{eff}}^k. $$
Finally, the capillary network is updated according to:}
$$\Lambda^{k}=\Lambda^{k-1}\cup \bigcup_{p\in \Lambda_d^k}[\xx_p^{k},\xx_p^{k-1}].$$
{If the criteria for branching are met two new capillary tips are added to $\Lambda_d^k$.}
	
\begin{rem}
{
{The update strategy} introduced here avoids the generation of capillary networks with a non-physiological structure, {such as segments shorter than a single cell length}. This also simplifies the meshing for the 1D domains, since it allows to obtain the 1D mesh as a simple sub-partition of the rectilinear segments composing the capillary network, avoiding the need of aggregation to ensure a lower bound to mesh size.}
\end{rem}

	\subsection{{Linearized time-discrete tumor cell fraction problem}}\label{subsec:discrete_tum_cell}

	System \eqref{eq:modelsummary_1} is {semi-}discretized in time on $\I_k$ using a one-step Backward Euler scheme, combined with an explicit treatment of the oxygen concentration, i.e.
	\begin{align}
	&\frac{\phic^k-\phi_c^{k-1}}{\Delta \I_k}=
\nabla \cdot (F_c(\phic^k)\nabla \phic^k)+S_c(\phic^k,c^{k-1})\phic^k  &\quad \text{ in } \Omega\label{eq:time_discrete_phic}\\
&\phi_l^k=\phimax-\phi_c^k&\quad \text{ in } \Omega\label{eq:time_discrete_phil}
	\end{align} 
		The resulting semi-implicit scheme avoids the implementation of a monolithic nonlinear 3D–1D coupled solver. Indeed, by evaluating $S_c$ in the known oxygen concentration $c^{k-1}$, the equation for $\phic^k$ becomes non linear only in $\phic^k$. This allows us to strongly reduce the computational complexity.
	 
	 Focusing on Equation \eqref{eq:time_discrete_phic} and concerning space discretization, let us denote by $\Phi$ and $\Phi^0$ respectively the discrete approximations of $\phic^k$ and $\phi_c^{k-1}$, such that
	$\Phi=\sum_{j=1}^N\Phi_j\eta_j$ and $\Phi^0~=~\sum_{j=1}^N\Phi^0_j\eta_j$ and let us introduce the matrices
	\begin{align*}
	\bm{K}(\Phi)\in \mathbb{R}^{N \times N} &\text{ s.t. }(K(\Phi))_{lj}=\int_\Omega F_c(\Phi)\nabla\eta_j\nabla\eta_l~d\bm{x}\\
	\bm{M}(\Phi)\in \mathbb{R}^{N \times N} &\text{ s.t. }(M(\Phi))_{lj}=\int_\Omega S_c(\Phi)\eta_j\eta_l~d\bm{x}\\
	\bm{B}\in \mathbb{R}^{N \times N} &\text{ s.t. }(B)_{lj}=\frac{1}{\Delta \I_k}\int_\Omega\eta_j\eta_l~d\bm{x}
	\end{align*}
	The fully discrete version of the equation describing the tumor cell fraction can be compactly written as
	\begin{equation}
G(\Phi)=0 ,\label{eq:non_linear}
	\end{equation} with
	$$G(\Phi):=(\bm{B}+\bm{K}(\Phi)-\bm{M}(\Phi))\Phi-\bm{B}\Phi^0.$$
	We solve numerically the non-linear equation \eqref{eq:non_linear} using Newton's method. At this aim we introduce the Jacobian matrix $\bm{J}\in \mathbb{R}^{N\times N}$ of $G(\Phi)$, whose elements are defined as
	\begin{align}
	J_{lm}(\Phi):=\frac{\partial G_l(\Phi)}{\partial \Phi_m}
	%&=\frac{\partial}{\partial \Phi_m}\Big[\sum_{j=1}^NB_{lj}(\Phi)\Phi_j+\sum_{j=1}^NA_{lj}(\Phi)\Phi_j-\sum_{j=1}^NM_{lj}(\Phi)\Phi_j\Big]=\nonumber\\
	%&=B_{lm}+\sum_{j=1}^N\frac{\partial A_{lj}(\Phi)}{\partial \Phi_m}\Phi_j+A_{lm}(\Phi)-\sum_{j=1}^N\frac{\partial M_{lj}(\Phi)}{\partial \Phi_m}\Phi_j-M_{lm}(\Phi)\nonumber\\
	&=B_{lm}+K_{lm}(\Phi)-M_{lm}(\Phi)+ (F_c'(\Phi)\nabla \Phi)\nabla \eta_l \eta_m-(S_c'(\Phi)\Phi)\eta_l\eta_m.
	\end{align}
	At each Newton iteration we then solve a linear problem in the form
	\begin{equation}
	\bm{J}(\Phi^{(n-1)})\delta\Phi=-G(\Phi^{(n-1)})
	\end{equation}
	with $\Phi^{(0)}=\Phi^0$,
	updating the vector of the tumor cell fraction degrees of freedom as
	$\Phi^{(n)}=\Phi^{(n-1)}+\delta \Phi$. Let us remark that the time interval $\mathcal{I}_k$ can be sub-partitioned into smaller sub-intervals if required for the convergence of the Newton scheme. 
    
\subsection{{Time-discrete 3D-1D pressure problem}}
	Problem \eqref{eq:modelsummary_2} involves the pressure of the liquid phase and includes a coupling between a 3D equation defined in the tissue and a 1D equation defined on the centerline of the capillary network. We tackle this coupling by means of the PDE constrained optimization-based strategy \cite{BGS3D1Ddisc, BGGPS2023}, as sketched in Section~\ref{sec:3D1D}.
	
	It can be noticed that $\phil^k$ can be computed from \eqref{eq:time_discrete_phic}-\eqref{eq:time_discrete_phil}, hence the first equation in \eqref{eq:modelsummary_2} is a quasi-static diffusion-reaction equation in $p$. We then rewrite it compactly as
	\begin{equation}
	\begin{cases}
		-\nabla \cdot \Big(\frac{\phil^k}{1-\phil^k}\frac{\kappa}{\mu}\nabla p^k\Big)+\phil^k\beta_p^{LS}\frac{S}{V}p^k=f_{3D}(\phic^k,\phil^k,c^{k-1})+\sum_{i \in Y^k}\limits \hat f_i^p(\hat p^k,\check{p}^k)\delta_{\Lambda_i^k} \quad \text{ in }\Omega\\
		-\frac{\pi R^4}{8 \mu}\frac{d^2\hat p_i^k}{ds^2}=-\hat f_i^p(\hat p^k,\check{p}^k)   	\hspace{5.8cm} \text{ in }(0,S_i^k) ~\forall i \in Y^k
	\end{cases}\label{eq:compact_p}
	\end{equation}
	with
	$$f_{3D}(\phic^k,\phil^k,c^{k-1})=\phil^k\beta_p^{LS}\frac{S}{V}p_{LS}-\phic^k S_c(\phic^k,c^{k-1})-\frac{1}{\Delta \I_k}(\phi_l^k-\phi_l^{k-1}).$$ 
Let us remark that the the oxygen concentration is taken at time $t_{k-1}$ to balance the reaction term in Equation \eqref{eq:modelsummary_1} and to decouple problems \eqref{eq:modelsummary_2} and \eqref{eq:modelsummary_3}.

The two equations in \eqref{eq:compact_p} can be decoupled by introducing proper auxiliary variables on the vessel segments $\Lambda_i^k$, namely $\psi_{\Omega,i}^k$ and $\psi_{\Lambda,i}^k$, approximating respectively the trace $\check p_i^k=\tr{p^k}{\Gamma_i^k}$ and the 1D pressure $\hat p_i^k$. We hence rewrite the problem as
\begin{equation}
\begin{cases}
-\nabla \cdot \Big(\frac{\phil^k}{1-\phil^k}\frac{\kappa}{\mu}\nabla p^k\Big)+\phil^k\beta_p^{LS}\frac{S}{V}p^k=f_{3D}(\phic^k,\phil^k,c^{k-1})+\sum_{i \in Y^k}\limits \hat f_i^p(\psi_{\Lambda}^k,\check{p}^k)\delta_{\Lambda_i^k} \quad \text{ in }\Omega\\
-\frac{\pi R^4}{8 \mu}\frac{d^2\hat p_i^k}{ds^2}=-\hat f_i^p(\hat p^k,\psi_{\Omega}^k)   	\hspace{5.8cm} \text{ in }(0,S_i^k) ~\forall i \in Y^k
\end{cases}\label{eq:opt_constraint}
\end{equation}
At each time step the auxiliary pressure variables $\psi_{\Omega}^k=\prod_{i \in Y^k}\psi_{\Omega,i}^k$ and $\psi_{\Lambda}^k=\prod_{i \in Y^k}\psi_{\Lambda,i}^k$ are chosen as the minimizers of the quadratic cost functional
$$J_p({\psi_{\Omega}^k},\psi_{\Lambda}^k)=\frac{1}{2}\sum_{i \in Y^k}\Big(||\check p_i^k-\psi_{\Omega,i}^k||^2_{\Lambda_i^k}+||\hat p_i^k-\psi_{\Lambda,i}^k||^2_{\Lambda_i^k}\Big)$$
constrained by \eqref{eq:opt_constraint}. As shown in \cite{BGS3D1Ddisc}, the corresponding discrete auxiliary pressures can be defined on two different sub-partions of the vessel segments, both independent from the partition used for the primary pressure variable $\hat{p}$. Here, for the sake of simplicity, we consider the same subdivision for both auxiliary pressures, discretizing them by linear finite elements as well.

To preserve the linearity of the constraint system, we use the auxiliary pressure variables at time $t_{k-1}$ to select the appropriate case when computing $\hat{f}_i^p$ (see \eqref{eq:fpi}).
The use of the auxiliary variables in this context is particularly advantageous when moving to the discrete problem, as we have decided to define the discrete auxiliary variables on the same subdivision of the capillary segments, and hence no interpolation is needed.

The algebraic formulation, as the one in equation~\eqref{eq:OPTDISCR}, can be obtained by assembling the integrals of the basis functions into matrices. The resulting discrete minimization problem is then solved at each time step by directly solving the associated system of optimality conditions.

\subsection{{Time-discrete 3D-1D oxygen problem}}
{Equation \eqref{eq:modelsummary_3} is semi-discretized in time by a backward Euler scheme, resulting in 
\begin{equation*}
\begin{cases}
\frac{\phi_l^{k}}{\Delta \I_k}(c^k-c^{k-1})-\nabla \cdot(\phil^k D_c\nabla c^k)+\bm{v}_l(p^k)\phil^k\cdot\nabla c^k\\\hspace{4.8cm}+ m_c \phil^k c^k(\phimax-\phil^k)=\sum_{i \in Y^k}\limits\hat f_i^c(\theta_{\Lambda}^k,\check{c}^k,)\delta_{\Lambda_i^k}\quad\text{ in } \Omega\\[1em]
			\frac{\pi R^2}{\Delta \I_k}(\hat{c}^k-\hat{c}^{k-1})-\pi R^2 \tilde{D}_c\frac{\partial^2\hat c^k }{\partial s^2}+\hat{v}_l(\hat{p}^k)\frac{\partial \hat c^k}{\partial s}=-\hat{f}_c^i(\hat{c}^k,\theta_{\Omega}^k) \hspace{3.3cm}	\text{ in }(0,S_i^k)~\forall i \in Y^k
			\end{cases}
			\end{equation*}
At each time step the auxiliary interface oxygen concentration variables $\theta_{\Omega}^k$ and $\theta_{\Lambda}^k$ are chosen as the minimizers of the cost functional
$$J_c(\theta_{\Omega}^k,\theta_{\Lambda}^k)=\frac{1}{2}\sum_{i \in Y^k}\Big(||\check c_i^k-\theta_{\Omega,i}^k||^2_{\Lambda_i^k}+||\hat c_i^k-\theta_{\Lambda,i}^k||^2_{\Lambda_i^k}\Big).$$ The discretization of $\theta_\Omega^k$ and $\theta_\Lambda^k$, as well as the assembly of $\hat{f}_c^i$ at the discrete level follow the same considerations discussed in the previous section for the auxiliary interface pressure variables and the assembly $\hat{f}_p^i$. }
\subsection{{Time-discrete VEGF problem}} \label{subsec:discrete_VEGF}
Equation \eqref{eq:modelsummary_4} is nothing but a 3D parabolic equation with singular sink term. We discretize it in time on $\I_k$ by a one time backward Euler scheme, obtaining
\begin{equation}
\frac{\phil^{k}}{\Delta \I_k}(g^k-g^{k-1})-\nabla \cdot(\phi_l^kD_g\nabla g^k)+\bm{v}_l(p^{k})\phil^k\cdot \nabla g^k+\sigma\phil^kg^k=\Gamma_g(\phic^k,c^k)-\sum_{i \in Y^k}\hat f_i^g(\check g^k)\delta_{\Lambda_i^k}\label{eq:time_discrete_g}.
\end{equation}
Let us observe how the adopted time discretization is fully implicit, as the value at time $t_k$ of all the other involved quantities is available from previous computations. Equation \eqref{eq:time_discrete_g} is discretized in space by linear finite element on the mesh $\mathcal{T}_h$.

	\section{Analysis of the sensitivity to parameters} \label{sec:results}

\begin{figure}
		\centering
		\includegraphics[width=0.4\linewidth]{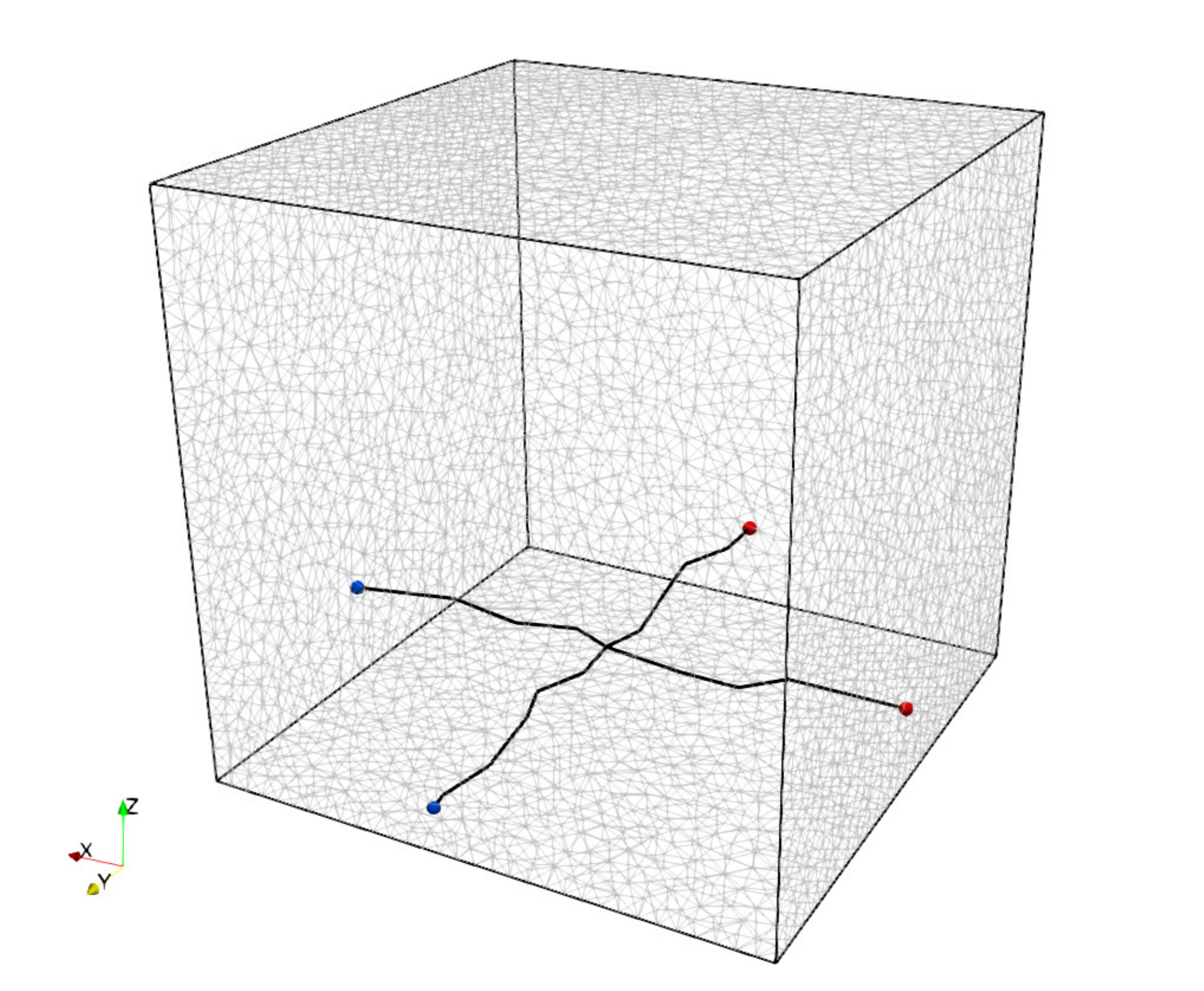}
		\caption{Initial capillary network $\Lambda^0$ and meshed computational domain $\Omega$.}
		\label{fig:initialgeometry}
	\end{figure}

	In this section we aim at analyzing the sensitivity of the proposed model to changes in the parameters. {We consider a representative three-dimensional test case consisting of a cubic tissue sample $\Omega$ in which a small initial capillary network is embedded (see Figure~\ref{fig:initialgeometry}). This configuration corresponds to the scenario labelled \textit{TestFace} in \cite{BGGPS2023}, where it was previously employed in the absence of the evolution of $\phi_c$ and $\phil$.}

	 Let $\bm{X}\in \mathbb{S}=[0,1]^K$ be a vector of input parameters, and let us denote by $Y(t,\bm{X})$ a given model output. 
     In order to analyse the sensitivity of $Y$ to variations in $\bm{X}$, we adopt the Elementary Effects (EE) method, originally introduced by Morris~\cite{Morris1991} and later refined by Campolongo et al.~\cite{Campolongo2007}. 
     %In order to analyze how sensitive $Y$ is to changes in $\bm{X}$ we adopt the Elementary Effects (EE) method presented by Morris in \cite{Morris1991}, in particular in the revised version proposed by Campolongo et al. in \cite{Campolongo2007}.
     Under this approach, we generate $R$ distinct trajectories in the input space $\mathbb{S}$.
	 Each trajectory consists of $K+1$ points, where the inputs are varied one at a time starting from a given initial point.
	 The region of experimentation is a $K$-dimensional $p$-level grid, meaning that each of the $K$ inputs can vary across $p$ different levels in the input space, with steps of size $\Delta$. Following \cite{Campolongo2007}, we choose $\Delta=\frac{p}{2(p-1)}$. 
	 
	 For a given point $\bm{X}\in \mathbb{S}$, the elementary effect associated to the $i$-th input is defined as 
	 $$d_{i}(t,\bm{X})=s_i\frac{Y(\bm{X}+s_i \bm{e}_i\Delta,t)-Y(\bm{X},t)}{\Delta}$$
	 where $s_i\in \lbrace -1,1\rbrace$, $\bm{e}_i=(e_i)_j\in \mathbb{R}^K$ such that $(e_i)_j=\delta_{ij}$ and $\bm{X}$ is randomly selected under the requirement that $\bm{X}+s_i\bm{e}_i\Delta\in \mathbb{S}$.
	 
	  In principle, the model could be run on all the $R$ trajectories in order to asses the analysis. Following \cite{Campolongo2007} we consider instead $r\ll R$ trajectories which are selected in order to maximize their dispersion in the parameter space (see \cite{Campolongo2007} for details). For each of the $r$ trajectories, we evaluate the model (i.e. we run a simulation) at all points and we compute the related elementary effects, ending up with $r$ independent elementary effects per input. In the more general case of non-uniform input distributions, trajectories are generated in the $[0,1]$ {uniform} probability space and then mapped back to the original distributions before evaluating the model. The step $\Delta$ always represents the sampling step in $[0,1]$ (see \cite{CampolongoCariboni2007} for examples).
	  
	   Finally, for each input we compute two sensitivity measures: the mean of the distribution of the absolute values of the elementary effects $\mu^*$ and the standard deviation of the elementary effects. The first measure evaluates the overall influence of the input on the output, with the absolute value preventing positive and negative effects from canceling out in case of non monotonic behavior (see \cite{Campolongo2007}). The second measure, instead, captures higher order effects, such as the interaction with other input parameters.

	 Tables \ref{table_geom}-\ref{table_growth} provide the set of input values employed in the simulations. For most of the parameters included in the sensitivity analysis, we assume uniform distributions over the reported intervals. Parameters that span several orders of magnitude are instead modeled using log-uniform distributions to reflect their scale-dependent variability (see Tables for details).
	 
	 \begin{table}\small
	 	\renewcommand*{\arraystretch}{1.25}
	 	\centering 
	 	\caption{Parameters defining the geometry}
	 	\label{table_geom}
	 	\begin{tabular}{|ccclc|}
	 		\hline
	 		\textbf{Parameter} &\textbf{Value} &\textbf{Unit}& \textbf{Description}&\textbf{Reference}\\
	 		\hline
	 		$L$ & $2.5$ &\small mm & Tissue sample  edge length& \cite{Gimbrone, angio_velocity, Cavallo} \\
	 		$R$ & $5\cdot 10^{-3}$ &\small mm &Vessel radius&\cite{capillary_wall}\\	\hline	
	 	\end{tabular}
	 \end{table}
	 \begin{table}\small 
	 	\renewcommand*{\arraystretch}{1.25}
	 	\centering 
	 	\caption{Parameters involved in the tumor cell problem. An interval $[a,b]$ indicates a uniform distribution between $a$ and $b$, while $\{a,b\}$ indicates a log-uniform distribution. 
	 	}
	 	\label{table_phic}
	 	\begin{tabular}{|ccclc|}
	 		\hline
	 		\textbf{Parameter} &\textbf{Value} &\textbf{Unit}& \textbf{Description}&\textbf{Reference}\\
	 		\hline
	 		$M$ & $1\cdot 10^{-4}$ & $\frac{\mathrm{mm}^2}{\mathrm{MPa} \cdot \mathrm{s}}$ &Motility parameter& \cite{GiversoGrillo}\\
	 		$E$ & $10$ &\small KPa & Young modulus& \cite{refE, GiversoGrillo}\\
	 		$\phimax$ & $1$ &\small $-$ &Maximum volume fraction of tumor cells& \cite{Byrne2003}\\
	 		$\phiz$ & $0.5$ &\small $-$ &Stress-free volume fraction of tumor cells& \cite{Preziosi, Ambrosi2003}\\
	 		$\gamma$ & $\bm{[1.16,~3.46]\cdot 10^{-2}}$ & $\frac{1}{h}$ \small &Tumor cell {proliferation} rate& \cite{Swanson2003, ByrneChaplain1996}\\	
	 		\multirow{1.7}{*}{$c_{\mathrm{ref}}$} & \multirow{1.8}{*}{$\bm{[8.5, 10.5]}$} & \multirow{1.7}{*}{\small mmHg} &{Threshold volumetric} oxygen concentration  &\\[-0.5em]&&&for tumor cell proliferation& \cite{Vaupel2004, Grimes}	\\\hline
	 	\end{tabular}
	 \end{table}
	 \begin{table}\small 
	 	\renewcommand*{\arraystretch}{1.25}
	 	\centering 
	 	\caption{Parameters related to the pressure problem. An interval $\{a,b\}$ indicates a log-uniform distribution between $a$ and $b$.}
	 	\label{table_pr}
	 	\begin{tabular}{|ccclc|}
	 		\hline
	 		\textbf{Parameter} &\textbf{Value} &\textbf{Unit}& \textbf{Description}&\textbf{Reference}\\
	 		\hline
	 		\multirow{1.5}{*}{$\beta_p^0$} & \multirow{1.5}{*}{$2.78\cdot 10^{-10} $} & \multirow{1.5}{*}{$\frac{\rm{mm}^2h}{\rm kg}$} &Hydraulic permeability &\\[-0.5em]&&&of healthy capillary wall&\cite{cattzun, Baxter_1}\\
	 		\multirow{1.5}{*}{$\Delta p_{onc}$} & \multirow{1.5}{*}{$25$ }& \multirow{1.5}{*}{\small mmHg} &Oncotic pressure jump &\\[-0.5em]&&& at the capillary wall&\cite{Levick, Levick_1, Jain1987}\\
	 		\multirow{1.5}{*}{$\beta_p^{LS}\frac{S}{V}$} &\multirow{1.5}{*}{$0.5$} & \multirow{1.5}{*}{\large $\rm \frac{mmHg}{h}$} &Effective permeability& \\[-0.5em]&&& of the lymphatic vessels&	\multirow{-1.5}{*}{\cite{cattzun,Baxter}} \\
	 		$\mu$ & $4 \cdot 10^{-3}$ & \small $\rm Pa \cdot s$& Blood viscosity&\cite{cattzun, Fung}\\
	 		$\tilde{p}_{in}$ & $33.75$& $\rm mmHg$ & inflow pressure&\cite{cattzun}\\
	 		$\tilde{p}_{out}$ & $35 $& $\rm mmHg$ & outflow pressure&\cite{cattzun} \\[0.5em]\multirow{1.5}{*}{$r_p^\beta$} & \multirow{1.5}{*}{$\bm{\lbrace 1,100\rbrace}$} & \multirow{1.5}{*}{$-$} & Permeability scale factor for&\\[-0.5em] &&& tumor-generated capillaries& \cite{Jain1987}\\
	 		\multirow{1.5}{*}{$\kappa$} & 	\multirow{1.5}{*}{$\bm{ \lbrace10^{-12} ,  10^{-7}\rbrace} $}& 	\multirow{1.5}{*}{\small$\rm mm^2$} &  Hydraulic permeability&\\[-0.5em]&&& of the tissue& \cite{Netti_Jain, GiversoGrillo}\\\hline
	 	\end{tabular}
	 \end{table}
	 \begin{table} \small 
	 	\renewcommand*{\arraystretch}{1.25}
	 	\centering 
	 	\caption{Parameters related to the oxygen problem. An interval $[a,b]$ indicates a uniform distribution between $a$ and $b$, while $\{a,b\}$ indicates a log-uniform distribution.}
	 	\label{table_oxy}
	 	\begin{tabular}{|ccclc|}
	 		\hline
	 		\textbf{Parameter} &\textbf{Value} &\textbf{Unit}& \textbf{Description}&\textbf{Reference }\\
	 		\hline
	 		\multirow{1.5}{*}{$\beta_c^0$} & \multirow{1.5}{*}{$12.6$} & \multirow{1.5}{*}{$\rm \frac{mm}{h}$} &Permeability of the healthy&\\[-0.5em] &&&capillary wall&\cite{cattzun, Baxter_1}\\
	 		$D_c$ & $4.86 $& $\rm \frac{mm^2}{h}$ &  Diffusivity, tissue &\cite{cattzun, Secomb1995}\\
	 		$\tilde{D}_c $ & $1.8\cdot 10^{3}$ & $\rm \frac{mm^2}{h}$ & Vascular diffusivity&\cite{cattzun, Grimes}\\
	 		$\tilde{c}_{in}$ & $95$& \small $\rm mmHg$ & inflow concentration& \cite{Levick, Levick_1, Grimes}\\[0.5em]\multirow{1.5}{*}{$r_c^\beta$} & \multirow{1.5}{*}{$\bm{\lbrace 1, 100\rbrace}$} & \multirow{1.5}{*}{$-$} & Permeability scale factor for&\\[-0.5em] &&& tumor-generated capillaries& \cite{Jain1987}\\
	 		\multirow{1.5}{*}{$m_c$} & \multirow{1.5}{*}{$\bm{[0.45 , 0.55]}$} &\multirow{1.5}{*}{$\rm \frac{1}{h}$} & Metabolization rate by&\\[-0.5em] &&&tumor cells& \cite{Grimes, Secomb1995}\\\hline
	 	\end{tabular}	
	 \end{table}
	 \begin{table}\small 
	 	\renewcommand*{\arraystretch}{1.25}
	 	\centering 
	 	\caption{Parameters related to the VEGF problem. Interval $[a,b]$ indicates a uniform distribution between $a$ and $b$.}
	 	\label{table_VEGF}
	 	\begin{tabular}{|ccclc|}
	 		\hline
	 		\textbf{Parameter} &\textbf{Value} &\textbf{Unit}& \textbf{Description}&\textbf{Reference}\\
	 		\hline
	 		$D_g$ & $0.18$ & $\rm\frac{{mm}^2}{h}$ &VEGF diffusivity& \cite{MacGabhann2006} \\
	 		$\sigma$ & $0.5$ & $\rm \frac{1}{h}$ &VEGF interstitial decay& \cite{Travasso} \\
	 		\multirow{1.5}{*}{$c^\star$} & \multirow{1.5}{*}{$11.5$} & \multirow{1.5}{*}{mmHg} & reference oxygen concentration&\\[-0.5em]&&& corresponding to $\Gamma_g =0.5G\phi_c$ & \cite{Vaupel2004, Grimes}\\
	 		$b$ & $11.5$ & \small - & dimensionless parameter in $\Gamma_g$& [-] \\\multirow{1.5}{*}{$\tilde{\sigma}$} &\multirow{1.5}{*}{$\bm{[0.2 , 2]}$} & \multirow{1.5}{*}{$\rm\frac{1}{h}$ } &Endothelial cell VEGF&\\[-0.5em]&&&  consumption rate & \cite{MacGabhann2006}\\
	 		$G$ &$\bm{[0.25 , 1]}$ & $\rm\frac{1}{h}$  &VEGF production rate& \cite{Travasso} \\\hline	
	 	\end{tabular}	
	 \end{table}
	 \begin{table}\small 
	 	\renewcommand*{\arraystretch}{1.25}
	 	\centering 
	 	\caption{Parameters related to capillary growth model. An interval $[a,b]$ indicates a uniform distribution between $a$ and $b$.}
	 	\label{table_growth}
	 	\begin{tabular}{|ccclc|}
	 		\hline
	 		\textbf{Parameter} &\textbf{Value} &\textbf{Unit}& \textbf{Description}&\textbf{Reference}\\
	 		\hline
	 		\multirow{2}{*}{$g_{lim}$} & 	\multirow{2}{*}{$0.25$}& 	\multirow{2}{*}{$10^{-13}\rm \frac{kg}{mm^3}$} &  minimum VEGF concentration & \\[-0.5em] &&& for endothelial cell proliferation & \multirow{-2}{*}{\cite{scirep_angionetwork,wang}}\\
	 		\multirow{1.5}{*}{$\bar{g}$} & 	\multirow{1.5}{*}{$\bm{[0.75 - 1.25]}$} & 	\multirow{1.5}{*}{$ 10^{-13}\rm \frac{kg}{mm^3}$} & VEGF concentration &{\cite{scirep_angionetwork,wang}}\\[-0.5em] &&& for $t_c=2\tau$& \\
	 		$l_e$ & $0.04$ & \small mm & endothelial cell length&\cite{sun}\\
	 		
	 		\multirow{1.5}{*}{$\alpha^w_{br}$} & \multirow{1.5}{*}{$\bm{[0.2-0.5]}$} & \multirow{1.5}{*}{-} & threshold of $\frac{||\bm{w}_\Pi||}{||\bm{w}||}$ &\\[-0.5em] &&&for branching &\cite{sun}\\
	 		$d_{br}$ & $\bm{[4.0,8.0]\cdot 10^{-2}}$ & \small mm & branching distance& \cite{andersonchaplain}\\
	 		\multirow{1.5}{*}{$\tau_{br}$} & \multirow{1.5}{*}{$\bm{[24,96]}$} & \multirow{1.5}{*}{\small h} & threshold age&\\[-0.5em] &&&for branching& \cite{andersonchaplain}\\		
	 		\multirow{1.5}{*}{$g_{br}$} & 	\multirow{1.5}{*}{$\bm{[\overline{g},2\overline{g}]}$} &	\multirow{1.5}{*}{$10^{-13}\rm \frac{kg}{mm^3}$} & VEGF concentration&\\[-0.5em]&&& for $P_{br}\approx1$ & [-]\\
	 	\multirow{1.5}{*}{$\tau$} & \multirow{1.5}{*}{$\bm{[12 , 48]}$} & \multirow{1.5}{*}{\small h} & endothelial cell proliferation &\\[-0.5em] &&&  parameter & \cite{SciannaPreziosi}\\ \hline
	 	\end{tabular}	
	 \end{table}
	 
	 Concerning the outputs, we monitor four time dependent quantities, namely \begin{itemize}
		\item {the relative tumor mass increase at time $t$,  $$P_\phi(t)=\frac{\mathcal{M}_{\phi_c}(t)-\mathcal{M}_{\phi_c}(0)}{\mathcal{M}_{\phi_c}(0)},$$
		where $\mathcal{M}_{\phi_c}(t)=\int_\Omega\phic(\xx,t)d\xx$;}
		\item the average quantity of oxygen available at time $t$
		$$C_\mathrm{avg}(t)=
		\frac{1}{|\Omega|}\int_{\Omega(t)}\phi_l(\bm{x},t)c(\bm{x},t)~d\bm{x}, $$
		\item the density $\rho_\mathrm{net}(t)$ of the capillary network, defined as
		$$\rho_\mathrm{net}(t)=\frac{L(t)}{N_{\rm{tips}}(t)},$$ where $L(t)$ denotes the total length of the capillary network of the capillary network at time $t$ and $N_{\rm{tips}}(t)$ is the number of capillary tips at time $t$;
		\item the fraction $V_\Omega(t)$ of vascularized tissue, which is approximated as the ratio between the number of elements intersected by the capillaries and the total number of elements in the mesh. 
	\end{itemize}
	
	All the simulations were performed on a tetrahedral mesh with $N=12438$ degrees of freedom.
    %, on a time interval 0 - 21 days, exporting the value of the outputs of interest at each day. 
    A uniform time discretization, with $\Delta\I_k=6 \rm{h}$ $\forall k=1,...,K$, was considered. At each time instant, the discrete primary 1D quantities ($\hat p$ and $\hat c$) were computed on the same partition of the capillary network, obtained by uniformly sub-partioning each capillary segment depending on the number of intersected tetrahedra. The four sets of auxiliary variables introduced to handle the 3D-1D coupled problems (two for the pressure and two for the oxygen) were also defined on a uniform grid, with half the number of nodes than the partition used for primary variables. Different choices are possible, as the discrete 3D-1D coupled problem is well-posed regardless the choice of the 1D partitions (see \cite{BGS3D1Ddisc}). {Two sensitivity analyses are performed. In SA1, we investigate the sensitivity of the model outputs to variations in nine parameters associated with all the equations of the model over a 21-day time span. In SA2, these nine parameters are fixed, and the analysis focuses instead on the parameters governing capillary geometry, again over a 21-day period. }

    {\subsection{SA1 test case}}
	\begin{figure}
		\centering
		
		\begin{subfigure}{0.9\linewidth}
			\centering
			\includegraphics[width=\linewidth]{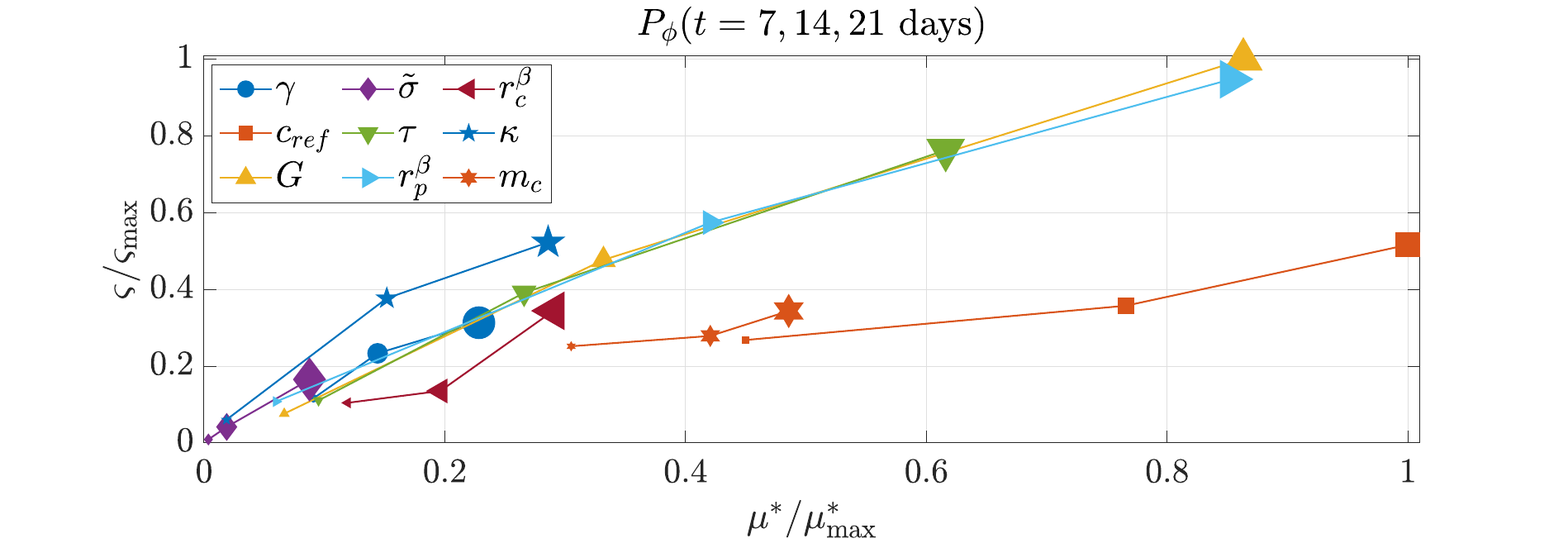}
			\caption{{Relative tumor mass increase.} Sensitivity measures normalized with respect to $\mu^*_\mathrm{max}=21.57$ and $\varsigma_\mathrm{max}=20.53$}
			\label{fig:pphi_7_14_21}
		\end{subfigure}
		
		\vspace{0.4cm} 
		
		\begin{subfigure}{0.9\linewidth}
			\centering
			\includegraphics[width=\linewidth]{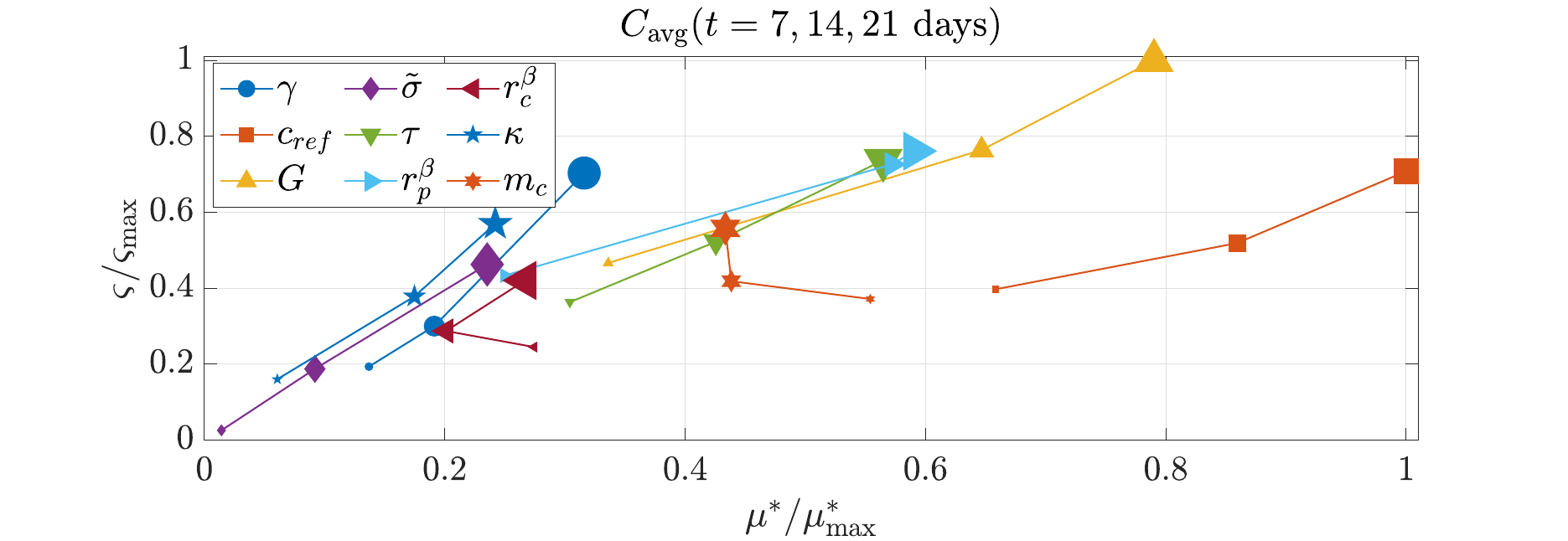}
			\caption{Average quantity oxygen. Sensitivity measures normalized with respect to $\mu^*_\mathrm{max}=1.41$ and $\varsigma_\mathrm{max}=1.23$}
			\label{fig:cavg_7_14_21}
		\end{subfigure}
		
		\caption{{SA1:} EE sensitivity tests at $t = 7, 14, 21$ days.}
		\label{fig:tum_quantities_7_14_21}
	\end{figure}
	
	\begin{figure}[ht]
		\centering
		\begin{subfigure}{0.9\linewidth}
			\centering
			\includegraphics[width=\linewidth]{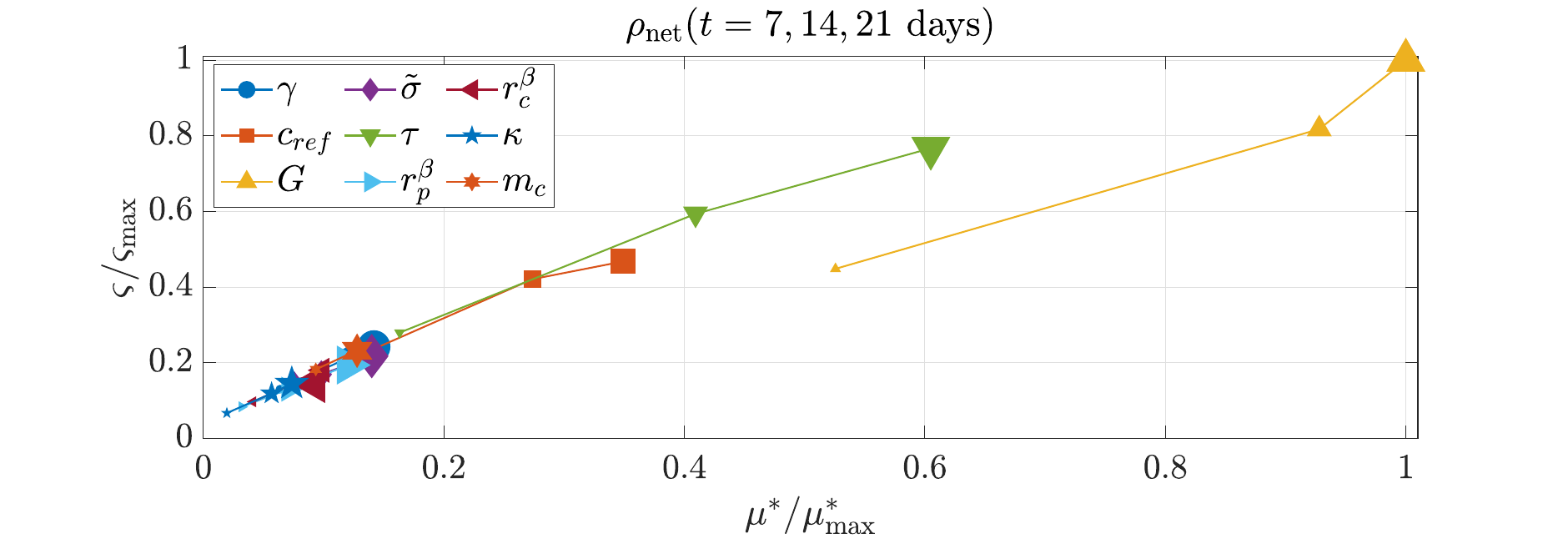}
			\caption{Density of the capillary network. Sensitivity measures normalized with respect to $\mu^*_\mathrm{max}=3.22$ and $\varsigma_\mathrm{max}=2.53$}
			\label{fig:rhonet_7_14_21}
		\end{subfigure}
		
		\vspace{0.4cm}

		\begin{subfigure}{0.9\linewidth}
			\centering
			\includegraphics[width=\linewidth]{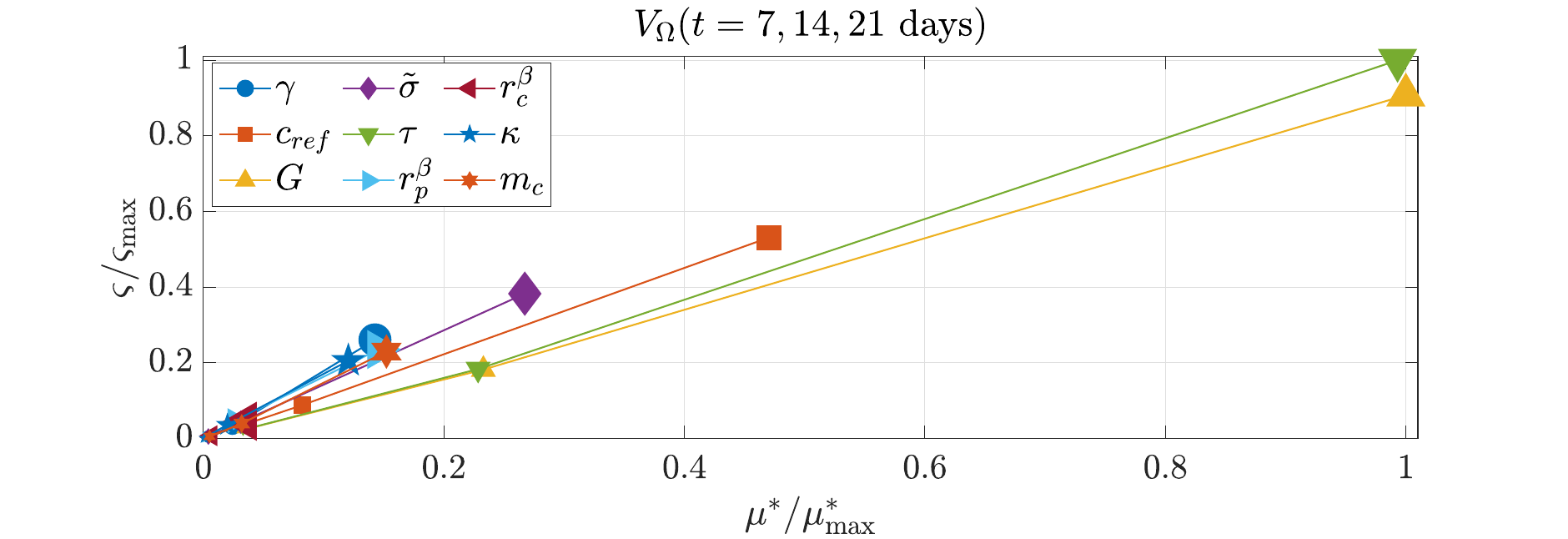}
			\caption{Fraction of vascularized tissue. Sensitivity measures normalized with respect to $\mu^*_\mathrm{max}=0.038$ and $\varsigma_\mathrm{max}=0.052$}
			\label{fig:vomega_7_14_21}
		\end{subfigure}
		
		\caption{{SA1:} EE sensitivity tests at $t = 7, 14, 21$ days.}
		\label{fig:net_quantities_7_14_21}
	\end{figure}
	
    {A first} sensitivity analysis is conduced on a 9-dimensional parameter space ($K=9$), allowing the inputs to vary across $p=4$ levels in the intervals reported in Tables \ref{table_phic}-\ref{table_growth}. {In particular we allow for the variation of $\gamma$, $c_{\mathrm{ref}}$ (Table \ref{table_phic}), $G$, $\tilde{\sigma}$ (Table \ref{table_VEGF}), $r_p^\beta$, $\kappa$ (Table \ref{table_pr}), $r_c^\beta$, $m_c$ (Table \ref{table_oxy}) and $\tau$ (Table \ref{table_growth}), whereas we fix the values of $\overline{g}=1\cdot 10^{13} \rm{kg/mm^3}$, $\alpha_{br}^w=0.3$, $d_{br}=4\cdot 10^{-2} \rm{mm}$, $\tau_{br}=48 \rm{h}$ and $g_{br}=1.5\cdot 10^{13} \rm{kg/mm^3}$.} 
    
    A total number of $R=1000$ trajectories was generated, from which $r=50$ trajectories were selected to evaluate the model. The results obtained at $t=7,14,21$ days are reported in Figures~\ref{fig:tum_quantities_7_14_21} and \ref{fig:net_quantities_7_14_21}, with the size of the markers being proportional to simulation time (i.e. small marker $\rightarrow$ $t=7$ days, big marker $\rightarrow$ $t=21$ days). For each output quantity, the values of $\mu^\star$ and $\varsigma$ are normalized by the maximum values obtained for that specific output, as the sensitivity indices are intended for qualitative comparison among parameters, and no cross-output comparison should be done using these figures.
	
	Figure \ref{fig:pphi_7_14_21} reports the results of the sensitivity tests for $P_\phi(t)$, i.e. the relative tumor mass increase. If we focus on the final time instant (big markers), we can observe a strong sensitivity to variations of parameters $G$, $r_p^\beta$, $\tau$ and $c_\mathrm{ref}$. These parameters are also the ones for which the sensitivity of $P_\phi(t)$ varies more over time. A strong sensitivity to $c_\mathrm{ref}$, already at an early simulation stage, is expected, as this parameter directly appears in the reaction term in Equation \eqref{eq:Sc}. The impact of $G$, $r_p^\beta$, $\tau$ is instead negligible at $t=7$ days, but their relevance (value of $\mu^*$) as well as their non-linear interaction with other parameters (value of $\varsigma$) increases rapidly with time. The sensitivity to $r_p^\beta$ is related to the fact that the permeability of tumor generated capillaries to fluids directly impacts on the fluid phase fraction and hence on the tumor cell fraction through \eqref{eq:algebraic}. The sensitivity to $G$, the rate at which tumor cells produce VEGF, and to $\tau$, which regulates the capillary growth rate, highlight how increased vascularization promotes tumor growth by enhancing tissue oxygenation, as it will be clear also from the comment to Figure \ref{fig:cavg_7_14_21}. It is interesting to note that $P_\phi(t)$ is marginally sensitive to {the tumor proliferation rate $\gamma$}. This behavior can be attributed to the fact that our analysis is restricted to fast-proliferating tumors, and therefore does not explore variations across multiple orders of magnitude of this parameter.
	
	Figure \ref{fig:cavg_7_14_21} shows the results obtained for the average oxygen concentration in the tissue, with similar considerations holding. Also in this case it is noteworthy to observe how oxygen availability is sensitive to parameters regulating capillary growth, highlighting the key role of vascularization in the overall process. It might be surprising that the average oxygen availability is more sensitive to $r_p^\beta$, than to $r_c^\beta$, which directly modulates the oxygen permeability of capillary walls. This behavior is explained by the fact that we are monitoring the average of $\phi_l c$, so that a variation in the liquid fraction has a particularly significant impact on the output. Similarly we explain the sensitivity to $c_\mathrm{ref}$, which impacts on $\phi_c$ and hence on $\phi_l$ through Equation \eqref{eq:algebraic}.
	
	The sensitivity results reported in Figure \ref{fig:net_quantities_7_14_21} refer instead to the outputs which are directly related to the capillary network development. In this case it is not surprising that both the network density $\rho_\mathrm{net}$ and the vascularized tissue fraction $V_\Omega$ are strongly influenced by $G$ and $\tau$. Since VEGF is produced by tumor cells, variations in $c_{\mathrm{ref}}$ also influence the formation and density of the capillary network. The influence of the VEGF metabolization rate $\tilde{\sigma}$ on $V_\Omega$ becomes influential towards the end of the simulation. Fast metabolization, indeed, induces a strong VEGF decrease around existing sprouts, promoting growth toward regions of higher VEGF concentration.
	
	To offer visual intuition on how parameter changes affect the simulation outcomes, we here propose some snapshots from two different runs in the set used for the sensitivity analysis (see table \ref{table_sets}).
    In Figure~\ref{fig:set1_set2} we compare the tumor cell fraction at $t=21$ days obtained with the parameters in Set~1 (Figure~\ref{fig:set1}) and Set~2 (Figure~\ref{fig:set2}){, as defined in Table \ref{table_sets}.} The simulation with Set~2 leads to a markedly higher proliferation, with an average tumor cell fraction of $\overline{\Phi}_c = 0.69$ and an extended region in which $\Phi_c$ exceeds this value. In contrast, the outcome of Set~1 yields a lower average of $\overline{\Phi}_c = 0.51$ and only a small region where the fraction is above this level. This is coherent with the fact that, in the Set 2 scenario, tumor cells proliferate at lower oxygen concentrations (see $c_{\mathrm{ref}}$ values), but also that the oxygen availability is in general higher. Indeed, in Figure \ref{fig:oxy_avg_set1set2} we can observe that, although the final average oxygen availability is very similar for the two parameter sets, the Set 2 scenario exhibits higher average availability during the whole simulation time. Concerning the morphology of the capillary network, we can observe already in Figure \ref{fig:set1_set2} how the results obtained in the Set 1 and in the Set 2 case look pretty similar. This is related to the fact that, although the rate $G$ at which tumor cells produce VEGF is lower in the Set 2 scenario, this is counterbalanced by the higher concentration of tumor cells observed in Figure \ref{fig:set2}. Despite being similar, the Set 2 capillary network ensures a better oxygen supply through the enhanced wall permeability (see factor $r_c^\beta$ and, in particular, factor $r_p^\beta$).\\

\begin{table}[h]\small
    \renewcommand*{\arraystretch}{1.3}
    \centering 
    \caption{{SA1:} Parameter values corresponding to the simulation of Figures \ref{fig:set1_set2}-\ref{fig:oxy_avg_set1set2}.}
    \label{table_sets_transposed1}
    \begin{tabular}{|l|ccccccccc|}
    \hline
         & $\gamma$ [$\text{h}^{-1}$] & $c_{\mathrm{ref}}$ [mmHg] & $G$ [$\rm{h}^{-1}$] & $\tilde \sigma$ [$\rm{h}^{-1}$] & $\tau$ [h] & $r_p^\beta$ [$-$] & $r_c^\beta$ [$-$] & $\kappa$ [$\rm mm^2$] & $m_c$ [$\rm h^{-1}$] \\
    \hline
    Set 1 & $1.93\cdot 10^{-2}$ & 10.5 & 1 & 1.4 & 36 & 21.38 & 1 & $3.22\cdot 10^{-9}$ & 0.55 \\
    Set 2 & $2.70\cdot 10^{-2}$ & 9.83 & 0.75 & 2 & 24 & 100 & 21.38 & $3.22\cdot 10^{-9}$ & 0.45\\
    \hline
    \end{tabular}
    \label{table_sets}
\end{table}
    
    \subsection{{SA2 test case}}
    {In this second analysis, the parameters in SA1 are fixed to the values reported in Set 1 (Table \ref{table_sets}), while we investigate the sensitivity of the outputs to variations in $K=5$ other parameters, namely $\overline{g}$, $\alpha_{br}^w$, $d_{br}$, $\tau_{br}$, and $g_{br}$ (Table \ref{table_growth}). These parameters regulate key aspects of the capillary network geometry, such as branching frequency, branching probability, and the spatial distance between branches. The sensitivity analysis is again carried out considering $R=1000$ total trajectories, from which $r=50$ were sampled to evaluate the model. Parameters were again allowed to vary across $p=4$ levels in their respective intervals of definition. 
     \begin{figure}[h]
        \centering
	\begin{subfigure}{0.5\linewidth}
		\centering
        \includegraphics[width=0.8\linewidth]{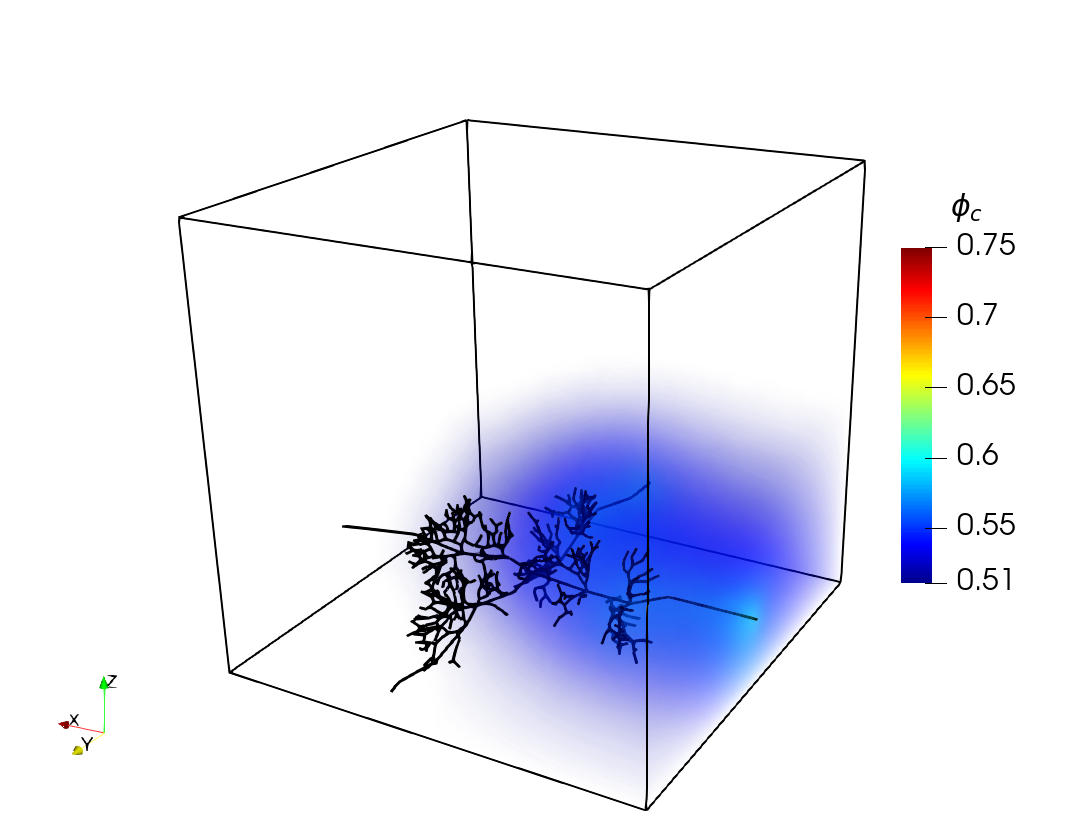}
        \caption{Parameters in set 1. $\overline\Phi_c=0.51$.}
        \label{fig:set1}
    \end{subfigure}%
    \begin{subfigure}{0.5\linewidth}
        \centering
\includegraphics[width=0.8\linewidth]{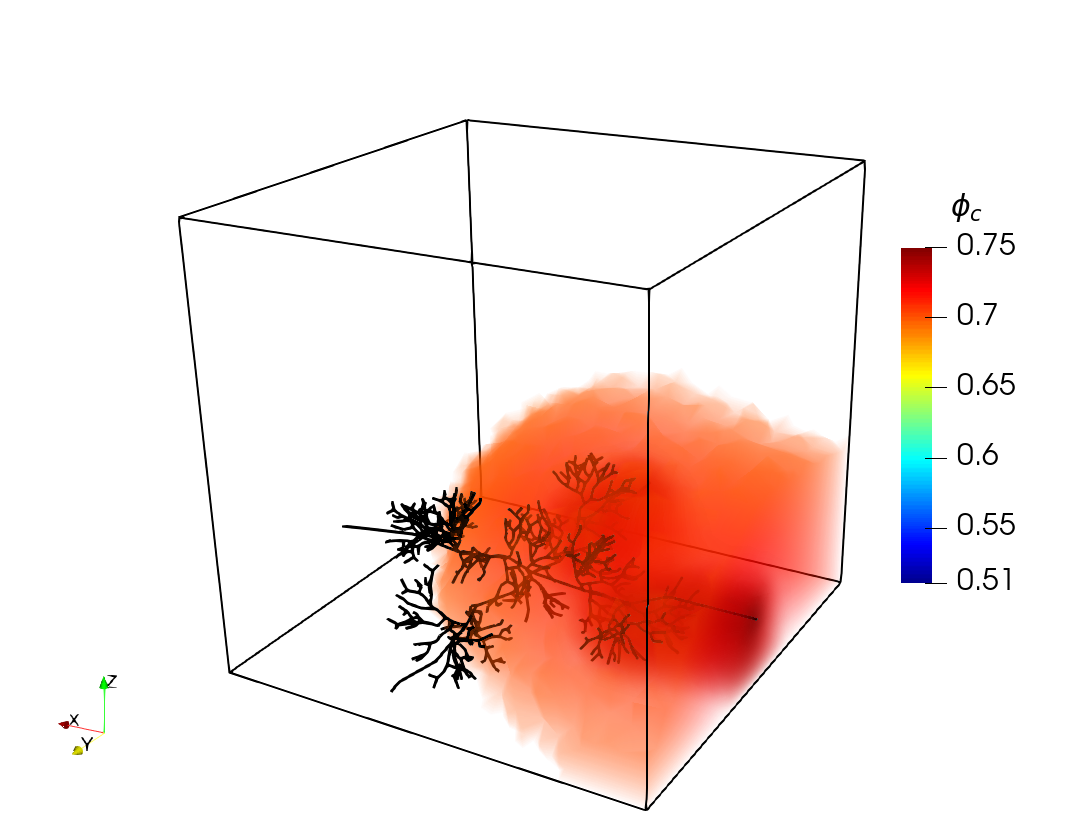}
		\caption{Parameters in set 2. $\overline\Phi_c=0.69$}
		\label{fig:set2}
        \end{subfigure}
        \caption{{SA1:} Tumor cell fraction at time $t=21$ days for two different parameter sets. White regions correspond to values of $\Phi_c$ below the average value $\overline \Phi_c$. }
        \label{fig:set1_set2}
	\end{figure}

    \begin{figure}[h]
        \centering
        \includegraphics[width=0.4\linewidth]{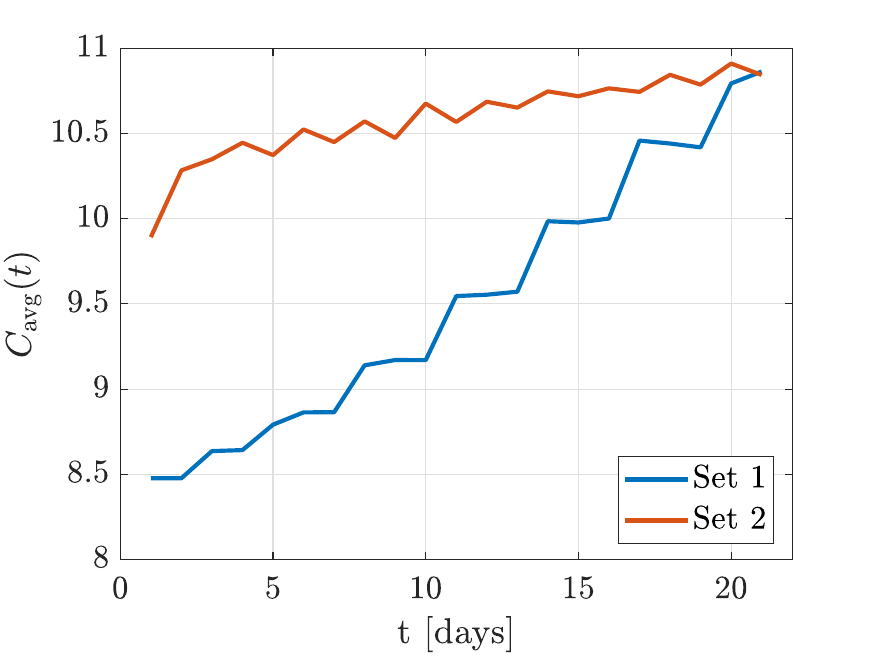}
        \caption{{SA1:}  oxygen concentration $C_{\mathrm{avg}} (t)$ with time for the Set 1 and Set 2 parameter sets.}
        \label{fig:oxy_avg_set1set2}
    \end{figure}
    
   Figures \ref{fig:SAgeom_cap_quantities}-\ref{fig:geom_tum_quantities_7_14_21} show the results of the sensitivity analysis for the different monitored outputs. Starting from Figure \ref{fig:SAgeom1_rhonet}, we observe, as expected, that the branching time $\tau_{br}$ strongly affects the network density. Interestingly, the impact of $\tau_{br}$ becomes apparent already at the early stages of the simulation. In the later stages, its direct effect remains roughly constant, while the interactions with the other parameters continue to grow, as reflected by the large value of $\varsigma$. This is consistent with the fact that branching frequency is only partially controlled by $\tau_{br}$: for example, even if the branching time is short, a low amount of VEGF may prevent branching from occurring.
The network density also strongly depends on the values of $\alpha_{br}^w$ and $d_{br}$. The former, being directly related to branching frequency, exhibits strong interactions with other parameters. In contrast, the values of $\varsigma$ for $d_{br}$ are lower, since this parameter regulates the distance between branches and interacts less with the others. Finally, $\overline{g}$ and $g_{br}$ have a comparatively smaller direct impact, but still exhibit notable interactions with the other parameters.
 Figure \ref{fig:SAgeom_Vomega} shows a strong sensitivity also of the fraction of vascularized tissue $V_\Omega(t)$ to both $\tau_{br}$ and $\alpha_{br}^w$. Notably, both the impact of these parameters and their nonlinear interactions with other factors increase with time.
 
    \begin{figure}[h]
        \centering
	\begin{subfigure}{0.48\linewidth}
		\centering
        \includegraphics[width=1\linewidth]{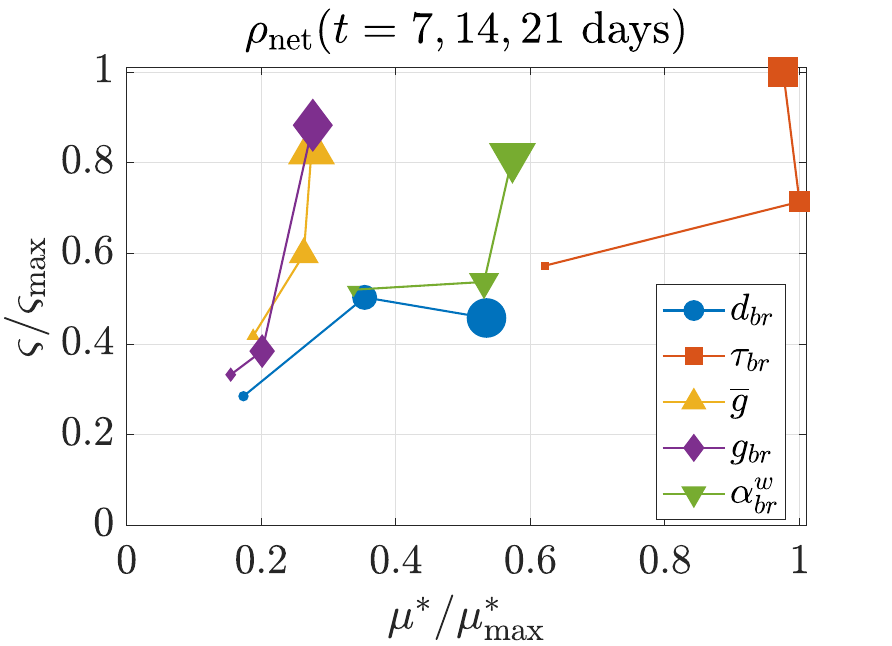}
        \caption{Density of the capillary network. Sensitivity measures normalized with respect to $\mu^*_\mathrm{max}=2.18$ and $\varsigma_\mathrm{max}=1.22$}
        \label{fig:SAgeom1_rhonet}
    \end{subfigure}\hfill%
    \begin{subfigure}{0.48\linewidth}
        \centering
		\includegraphics[width=1\linewidth]{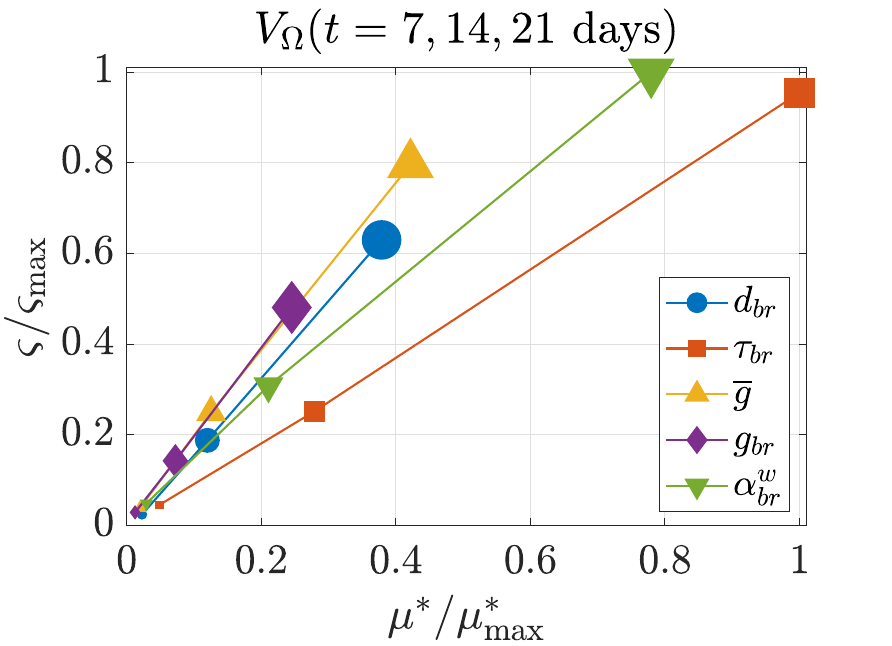}
		\caption{Fraction of vascularized tissue. Sensitivity measures normalized with respect to $\mu^*_\mathrm{max}=0.019$ and $\varsigma_\mathrm{max}=0.013$}
		\label{fig:SAgeom_Vomega}
        \end{subfigure}
          \caption{{SA2:} EE sensitivity tests at t=7,14,21 days.}
        \label{fig:SAgeom_cap_quantities}
	\end{figure}
    \begin{figure}[H]
        \centering
	\begin{subfigure}{0.48\linewidth}
		\centering
\includegraphics[width=1\linewidth]{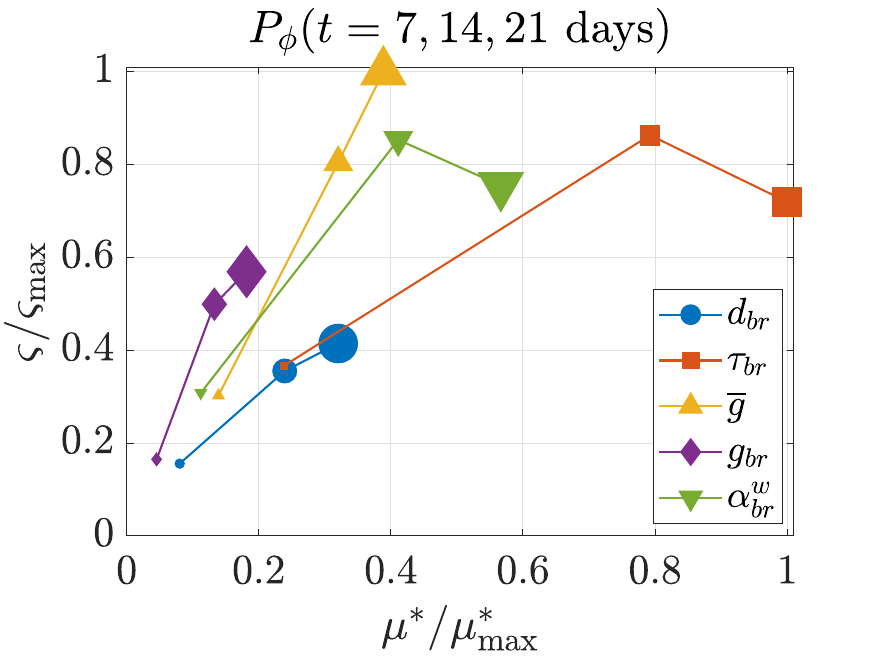}
        \caption{Relative tumor cell fraction increase. Sensitivity measures normalized with respect to $\mu^*_{\rm{max}}=11.50$ and $\varsigma_{\rm{max}}=4.26$}
        \label{fig:geom_pphi_7_14_21}
    \end{subfigure}\hfill%
    \begin{subfigure}{0.48\linewidth}
        \centering
\includegraphics[width=1\linewidth]{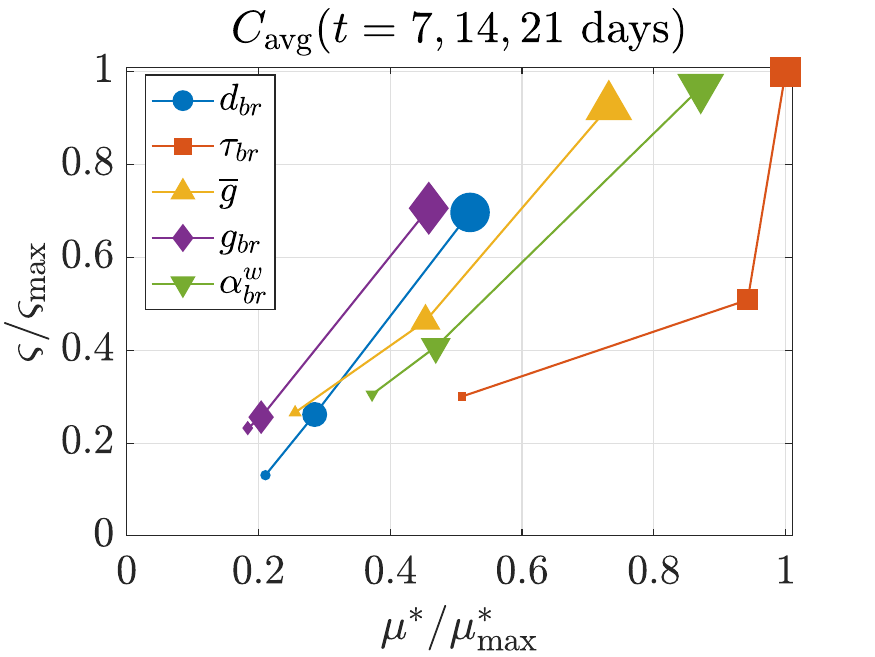}
		\caption{Average quantity oxygen. Sensitivity measures normalized with respect to $\mu^*_\mathrm{max}=0.46$ and $\varsigma_\mathrm{max}=0.54$}
		\label{fig:geom_cavg_7_14_21}
        \end{subfigure}
        \caption{SA2: EE sensitivity tests at t=7,14,21 days.}
        \label{fig:geom_tum_quantities_7_14_21}
        \end{figure}

Concerning tumor mass growth, the sensitivity analysis results are shown in Figure \ref{fig:geom_pphi_7_14_21}, with the branching time $\tau_{br}$ being again the most influential parameter. To better highlight the impact of the capillary geometry on tumor proliferation, Figure \ref{fig:geom:phicintime} shows the tumor cell volume fraction over time for two representative runs from the sensitivity analysis dataset (see Table \ref{table_sets1ab}). The corresponding final capillary networks are shown in Figure \ref{fig:set1a_set1b}, revealing one network that is considerably denser and another that is much less developed. In the more developed network, the final average tumor cell fraction is approximately 
%nine percent 
$9 \%$ higher (see Figure \ref{fig:geom:phicintime}), which represents a substantial difference considering that all other parameters are held constant. The impact of the capillary geometry on tumor growth is again related to the increase of oxygen availability. The outcomes of the sensitivity analysis for $C_{avg}$ are shown in Figure \ref{fig:geom_cavg_7_14_21}. Also in this case, it is noteworthy that the influence of $\tau_{br}$ increases rapidly during the early stages of the simulation, while in the later stages its direct effect remains almost constant and the interactions with the other parameters become more evident.}

\begin{table}[h]\small
    \renewcommand*{\arraystretch}{1.2}
    \centering 
    \caption{{SA2:} Parameter values corresponding to the simulation of Figures \ref{fig:set1a_set1b}-\ref{fig:geom:phicintime}.}
    \begin{tabular}{|l|ccccc|}
    \hline
         & $d_{br}$ [mm] & $\tau_{br}$ [h] &$\overline{g}$ [$10^{-13}\rm{kg/mm^3}$] & $g_{br}$ [$10^{-13}\rm{kg/mm^3}$] & $\alpha_{br}^w$ [-]    \\
    \hline
    Set 1.a & $6.67\cdot 10^{-2}$ & 48 & 1.08 & 1.81 & 0.3 \\
    Set 1.b & $6.67\cdot 10^{-2}$ & 96 & 1.08 & 1.81 & 0.5 \\
    \hline
    \end{tabular}
    \label{table_sets1ab}
\end{table}

 \begin{figure}[h]
        \centering
	\begin{subfigure}{0.45\linewidth}
		\centering
        \includegraphics[width=1\linewidth]{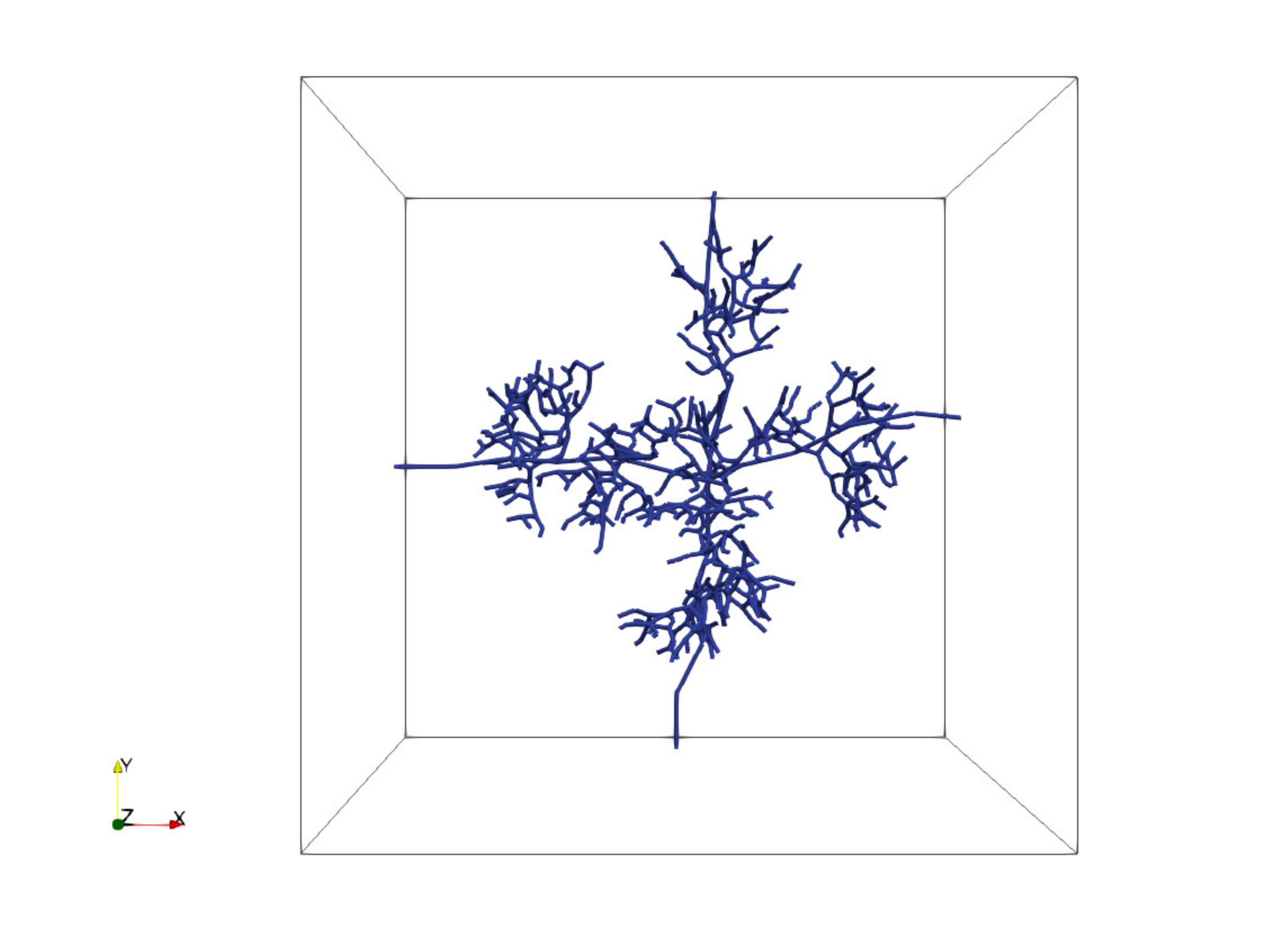}
        \caption{Parameters in set 1.a}
        \label{fig:set1a}
    \end{subfigure}%
    \begin{subfigure}{0.45\linewidth}
        \centering
		\includegraphics[width=1\linewidth]{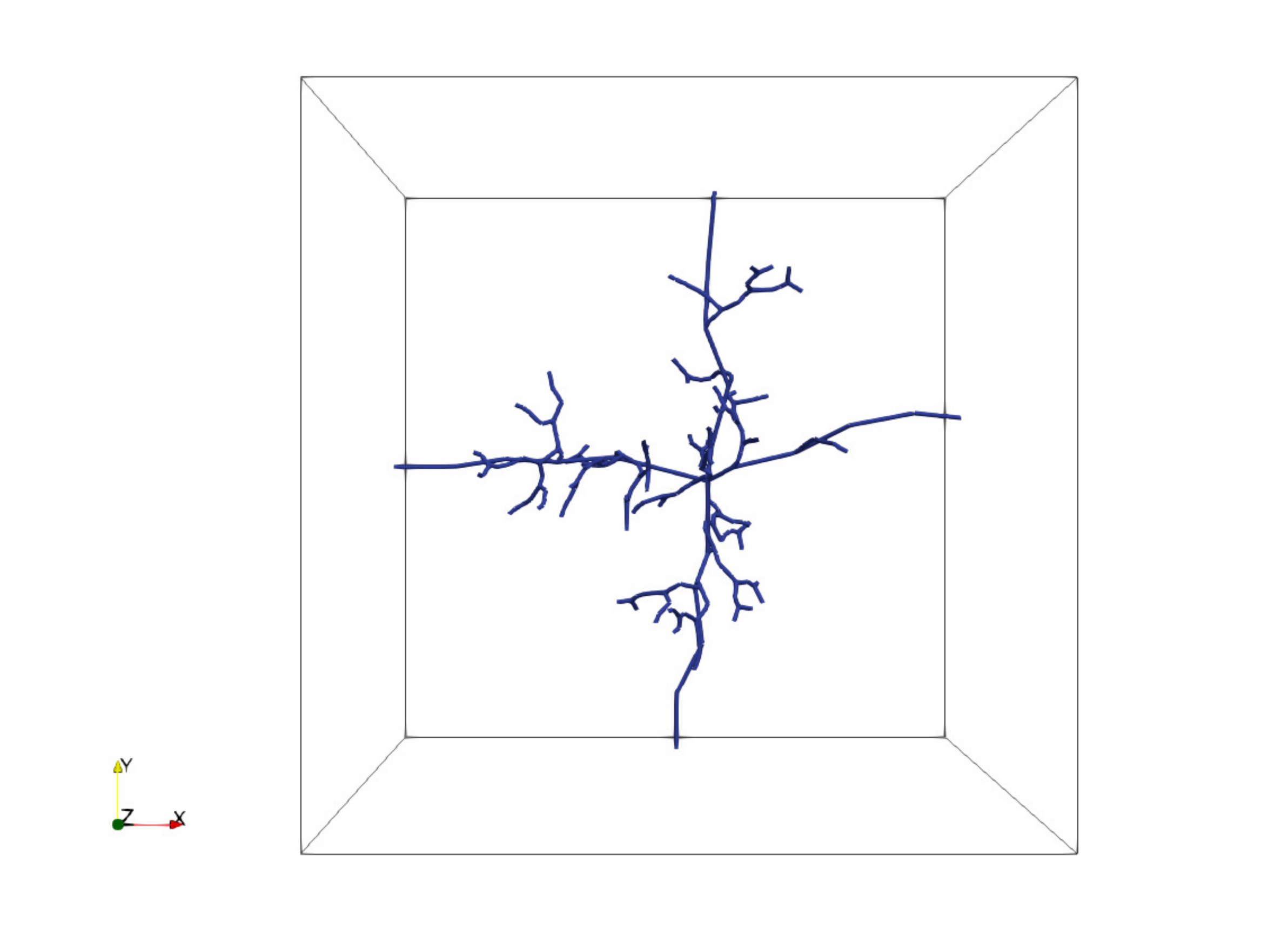}
		\caption{Parameters in set 1.b}
		\label{fig:set1b}
        \end{subfigure}
        \caption{SA2: Capillary network configuration at time $t=21$ days for two different parameter sets, {as defined in Table \ref{table_sets1ab}}.}
        \label{fig:set1a_set1b}
	\end{figure}
\begin{figure}[h]
    \centering
    \includegraphics[width=0.45\linewidth]{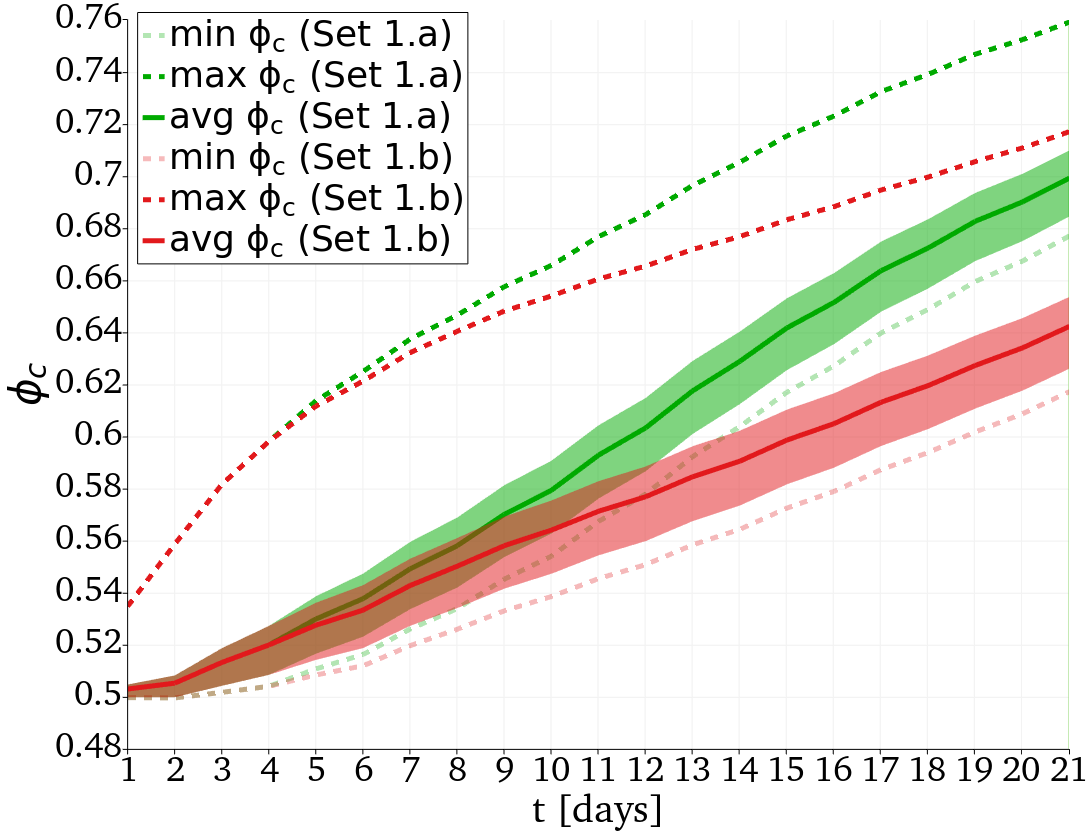}
    \caption{SA2: Tumor cell volume fraction statistics for two different parameter sets, as defined in Table \ref{table_sets1ab}.}
    \label{fig:geom:phicintime}
\end{figure}

\section{Conclusions} \label{sec:conclusions}
{This work extends our previous modeling framework for tumor-induced angiogenesis \cite{BGGPS2023} by incorporating a multiphase description of the tumor tissue. In the proposed model, tumor cell proliferation is related to local oxygen concentration, which increases in regions perfused by the growing vascular network, resulting in a nonlinear feedback between cell volume fraction, oxygen distribution, and vascular dynamics.

From a methodological standpoint, the paper exploits an optimization-based 3D–1D coupling strategy to model the exchanges of fluids and oxygen through the capillary wall. This coupling approach, in which the problem is formulated as a constrained minimization problem, leads to a well-posed discrete problem, enabling the use of independent and non-conforming discretizations for all the different variables. This aspect represents a  key advantage of the approach, particularly in the presence of evolving geometries generated by the hybrid tip-tracking algorithm, which dynamically updates the vascular architecture in response to angiogenic stimuli.

The simulations highlight how the evolution of the capillary network strongly affects tissue perfusion and oxygen availability, which in turn regulate tumor dynamics through nonlinear feedback mechanisms. The performed sensitivity analysis, based on the Morris elementary effects method, provides further insight into the relative importance of model parameters, confirming the central role of vascular-related quantities and supporting the use of the model as a tool for parameter investigation and calibration.

Several directions for future research can be envisaged. First, the framework will be extended to include the simulation of metastasis by allowing tumor cells to penetrate the capillary wall. This will be achieved by extending the 3D–1D coupling framework to the equation describing tumor cell dynamics, accounting for cell intravasation processes. Second, the introduction of vascular remodeling mechanisms will enable the study of structural changes in the capillary network, with particular relevance for modeling and assessing the effects of therapeutic interventions. Finally, the model will be further enriched by introducing a distinct interface between the tumor and surrounding healthy tissue, allowing the explicit simulation of interactions across this boundary. This extension will enable a more realistic representation of the tumor microenvironment and provide a stronger basis for comprehensive simulations.}

\section*{Acknowledgments}
{The authors wish to thank Annachiara Colombi (DISMA, Politecnico di Torino) for her kind and valuable advice.}
Author Denise Grappein acknowledges that part of this work was conducted while holding a postdoctoral fellowship funded by the Italian Institute of High Mathematics (INdAM), hosted at the research unit of Politecnico di Torino.

	\printbibliography
\end{document}